\documentclass[10pt,a4paper]{amsart}
\usepackage{amsmath}
\usepackage{rotating}
\def\hh{\mathcal{H}}
\def\s3{\mathcal{S}_3}
\def\dimo{\textit{Proof}$\quad$}

\def\ss{\star}

\def\ot{{\oplus 3}}
\def\os{{\oplus 6}}

\def\Torus{\mathbb{T}}

\def\me{\mathbf{g}}
\def\contr{\neg}
\def\dime{3}

\def\contr{\rightharpoonup}

\def\R{\mathbb{R}}
\def\C{\mathbb{C}}
\def\Z{\mathbb{Z}}

\def\gg{\mathcal{L}_{\mathbf{3},\C}}
\def\ggr{\mathcal{L}_\mathbf{3}}
\def\so3{\mathbf{so}(3,\R)}
\def\sl3c{\mathbf{sl}(3,\C)}

\def\ka{K\"{a}hler\;}

\newtheorem{teo}{Theorem}[section]

\newtheorem{cor}[teo]{Corollary}

\newtheorem{lem}[teo]{Lemma}
\newtheorem{pro}[teo]{Proposition}
\newtheorem{dfn}[teo]{Definition}

\newtheorem{rmk}[teo]{Remark}

\title{A natural Lie superalgebra bundle on rank three WSD manifolds}
\author{Giovanni Gaiffi, ~Michele Grassi}
\date{July 11, 2007}
\begin{document}

\begin{abstract}
We determine the structure of the  $*$-Lie superalgebra generated by a set of carefully chosen natural operators of an orientable WSD manifold of rank three. This Lie superalgebra is formed by global sections of a  natural Lie superalgebra bundle, and turns out to be a product of $\mathbf{sl}(4,\C)$ with the full special linear superalgebras of some graded vector spaces isotypical with respect to a natural action of $\mathbf{so}(3,\R)$.  We provide an explicit description of one of the real forms of this superalgebra, which is geometrically natural being  made of   $\mathbf{so}(3,\R)$-invariant operators which preserve the Poincar\'e (odd Hermitean) inner product on the bundle of forms. 
\end{abstract}
\maketitle
\section{Introduction}
In the present paper we continue the search for natural ($*$-)Lie superalgebra bundles on Weakly Self Dual manifolds, started in \cite{GG}. By carefully choosing the initial geometrical operators canonically defined on any orientable Weakly Self Dual manifold of rank three, we generate a  complex $*$-Lie superalgebra $\gg$ which is invariant with respect to pointwise automorphism group $\mathbf{Aut}_p(X) = \mathbf{so}(3,\R)$ at any point $p$ of the WSD structure. We describe the action of this $*$-Lie superalgebra on the differential forms of the WSD manifold, decomposing them into irreducible representations, and we also show that the algebra of all the  operators of $\gg$  which preserve the Poincar\'e (odd Hermitean) inner product on the bundle of forms at any point $p$ is a real form of it. \\
We prove the following main theorem:\\

\textbf{Theorem~\ref{teo:mainteo}} (simplified form)\textit{On any rank three orientable WSD manifold there is a natural real Lie superalgebra bundle, based on the real Lie superalgebra
\[\mathbf{su}(20|20,<~,~>)~\oplus~ \mathbf{su}(36|36,<~,~>)~\oplus~ \mathbf{su}(20|20,<~,~>)~\oplus~ \mathbf{sl}(4,\R)\]}

\noindent where we indicate with $\mathbf{su}(n|n,<~,~>)$ the sub Lie superalgebra of $\mathbf{sl}(n|n,\C)$  which preserves the standard odd Hermitean inner product $<~,~>$ on $\C^{n|n}$. This is not the "usual" $\mathbf{su}(n|n)$, as the inner product is odd; indeed, in a good basis the matrices of this Lie superalgebra have the form
\[\begin{pmatrix} A & B =~ ^t{\overline{B}}\\ C =~ -{^t{\overline{C}}} & -{^t{\overline{A}}}\end{pmatrix}\]
The motivations for this search go back to \cite{G2}, where it was shown that on any weakly self dual manifold  (WSD for brevity) there is a natural $\mathbf{sl}(4,\R)$ bundle. Afterwards in \cite{G3} it was conjectured that there should be a "good" (super)Lie algebra bundle acting on $\bigwedge^*T^*X$ for a rank three WSD manifold $X$ (in that case degenerate of dimension $11$), which should give rise, using the (generalized) lagrangian dynamics of \cite{G1}, to a natural and geometrically motivated field theory on any such WSD maniflold. It was further conjectured that once quantized this field theory would provide a good playground where to look for a unifying theory for the various string theories. To put this last part of the conjectured construction into context, one should recall that the generalized lagrangian dynamics of \cite{G1} is also as well a generalized hamiltonian dynamics, at least for nondegenerate lagrangians, and this is clearly the starting point in the quest for a quantization procedure. \\
Our basic motivations come therefore from mathematical physics, and from this point of view it is both reassuring and stimulating that our final result Theorem \ref{teo:mainteo} finds a Lie superalgebra of "super-antihermitean" operators with respect to an odd nondegenerate Hermitean inner product $<~,~>$. \\
\noindent Polysymplectic manifolds were introduced in \cite{G1} to provide a geometric approach to string theory, and WSD manifolds, which are very special polysymplectic manifolds,  were introduced in \cite{G2} to build a geometric approach to Mirror Symmetry. In just a few words, if two Calabi-Yau manifolds are mirror dual, it is conjectured that there should be a family of WSD manifolds which converges (in a normalized Gromov-Hausdorff sense) to either of them on two different boundary points of the deformation space, and that for this to happen for any family of WSD manifolds it is actually necessary to have a mirror pair of Calabi-Yau manilfolds  (see \cite{G2},\cite{G3} and the introduction to \cite{GG} for further details).  In \cite{G2} this approach was shown to work for simpler cases, like Elliptic curves and Affine-K\"ahler manifolds, while in \cite{G3} it was shown that many aspects of this method do work in the case of anticanonical families in complex projective spaces, which includes the historically significant case of the quintic threefold. \\
Moved by these motivations, in \cite{GG} we started with a "proof of concept" computation for rank two WSD manifolds, to develop the tecniques necessary to tackle the higher rank cases. Even if  the rank two context was clearly less significant from a physical point of view, the result was nevertheless interesting algebraically and geometrically, as we found a natural $\mathbf{sl}(6,\C)$ bundle, which is a direct generalization of the $\mathbf{sl}(2,\C)$ bunlde arising naturally on K\"ahler manifolds.\\
As mentioned above, here we show that indeed  in the rank three nondegenerate case (and therefore also in the degenerate one) there is a very structured Lie superalgebra bundle acting as conjectured. We do not provide suggestions for the possible generalized lagrangian theory based on this bundle, even if there are some natural candidates which one should try first, from a geometrical point of view. \\
To arrive to the final description, we faced a much more intricated situation with respect to \cite{GG}, as first of all in the rank three (oriented) case  one is working with superLie algebras, instead of ordinary (even) Lie algebras. To see why this happens, it is necessary to give the definition of WSD manifold: we reproduce here the rank three nondegenerate case, which is the one relevant for our present purposes (look at \cite{G2} for the general definition); also recall that we always assume the WSD manifolds to be orientable:
\begin{dfn}
A {\em rank three nondegenerate weakly self-dual manifold} (WSD manifold for brevity) is
given by a smooth manifold $X$, together with two smooth $2$-forms
$\omega_1,\omega_2$ a Riemannian metric and a third smooth
$2$-form $\omega_D$ (the {\em dualizing} form) on it, which satisfy the following conditions:\\
1) $d\omega_1 = d\omega_2 = d\omega_D = 0$ and the distribution $\omega_1^0+\omega_2^0$ is integrable.\\
2) For all $p\in X$ there exist an orthonormal basis
\[dx_1,dx_2,dx_3,dy^1_1,dy^1_2,dy^1_3,dy^2_1,dy^2_2,dy^2_3\]
 of  $T_p^*X$ such that the
\[(\omega_1)_p = \sum_{i=1}^m dx_i\wedge dy^1_i,~ ~(\omega_2)_p =
\sum_{i=1}^m dx_i\wedge dy^2_i,~ ~(\omega_D)_p = \sum_{i=1}^m
dy^1_i\wedge dy^2_i\] 
Any 
basis of $T_pX$ dual to a basis  of $1$-forms  as above  is said to
be {\em adapted} to the structure, or {\em standard}
\end{dfn}
\noindent
Starting from this geometric data, one constructus six natural "wedge" operators
\begin{dfn}\label{dfn:Loperators}
For $\phi\in
\Omega^*_\C X$,
\[L_0(\phi) = \omega_D\wedge \phi,\qquad L_1(\phi) = -\omega_2\wedge
\phi,\qquad L_2(\phi) = \omega_1\wedge \phi\] 
{\small \[V_0(\phi) = dx_0\wedge dx_1\wedge dx_2\wedge \phi,\qquad V_1(\phi) = dy^2_0\wedge dy^2_1\wedge dy^2_2\wedge \phi, \qquad V_2(\phi) = dy^1_0\wedge dy^1_1\wedge dy^1_2\wedge \phi\]}
\end{dfn}
\noindent
and their (super) adjoints $\Lambda_j = L_j^\ss$, $A_j = V_j^\ss$ with respect to the standard Hermitean form $(~,~)$. Here the $V_j$ are actually wedge operators with volume forms associated to canonical oriented distributions on $T^*X$, and therefore are well defined. As can be seen by inspection, while the $L_j$ are even, the $V_j$ are odd, while in the rank two case the corresponding $V_j$ operators were also even. We call $\gg$ the complex Lie superalgebra generated by these $12$ operators. One of the nice aspects of this theory is that  the operators above are invariant with respect to a group which is the space of sections of the $\mathbf{SO}(3,\R)$ bundle determined point by point by the automorphism group $\mathbf{Aut}_p(X)$ of the WSD structure. The Lie algebra associated to this group acts on the faithful, defining representation of $\gg$, namely $\bigwedge^*T^*X$, and by conjugation on $\gg$. This action on both $\gg$ (which as mentioned is $\mathbf{SO}(3,\R)$-invariant) and on its representation makes the representation theory that we use both richer and  more complex. This is a general trend, and in the rank $k$ case the invariance group is the space of sections of a $\mathbf{SO}(k,\R)$ bundle. To obtain the complete description given in Theorem \ref{teo:mainteo}, in the following corollary and in the tables in appendix, we used also the $\mathbf{sl}(4,\C)$ structure uncovered in \cite{G1}, which is present on a WSD manifold of any rank. To obtain enough operators we needed however to use also the supercommutators $[V_j,A_j]$ to split the operators in this $\mathbf{sl}(4,\C)$ into homogeneous components.\\
Another geometric ingredient that comes into play is the odd, superHermitean, nondegenerate Poincar\'e inner product $<~,~>$ between forms, given at any fixed point $p\in X$ by the formula
\[<\alpha_p,\beta_p>dVol_p = \alpha_p\wedge\overline{\beta_p}\] 
where $dVol_p$ is the volume form at the point $p$ given by the metric and the orientation. 
We show that the the operators $\imath L_j,\imath\Lambda_k,\imath V_s, A_t$ are all superantiHermitean with respect to the inner product $<~,~>$, and  the real Lie superalgebra $\ggr$ which they generate    is a proper real  form for $\gg$.
Moreover, the inclusion of $\ggr$ inside the full $\mathbf{SO}(3,\R)$-invariant superantiHermitean superalgebra (which has dimension 8444) is actually almost an isomorphism, its dimension being off by a mere $48$. \\

Coming to a more precise description of the contents of this paper, in Section \ref{sec:alginv} we introduce the algebra and its symmetries. There are many pieces of the puzzle that need to be introduced, and which will be put together in the subsequent sections. Among them, besides the geometric operators $L_j,\Lambda_j,V_j,A_j$ and the real form $\ggr$, we have the operators $J_k$ generating the $\mathbf{SO}(3,\R)$ symmetry, the twisted adjoint $\phi^*$ to an operator, and the Poincar\'e odd Hermitean inner product $<~,~>$, which we show that is preserved by $\ggr$. There is also the Hodge star, which is ubiquitous and sends even forms to odd ones, as our manifolds are odd dimensional.\\
In Section \ref{sec:so2R}, after introducing a natural action of the permutation group $\mathcal{S}_3$, we first prove that $\ggr$ commutes with $\mathbf{SO}(3,\R)$ and then we continue by decomposing (doing a form of plethysm for $\mathbf{so}(3,\R)\otimes\C \cong \mathbf{sl}(2,\C)$) the representation $\bigwedge^*_\C T^*X$ with respect to the action of $\mathbf{sl}(2,\C)$. Any isotypical $\mathbf{sl}(2,\C)$-module will be clearly an invariant $\gg$-submodule from Schur's lemma; even more efficiently, the highest weight vectors inside an isotypical component form an  invariant $\gg$-submodule. In Corollary \ref{cor:subend} we make this precise, by mapping $\ggr$ inside a direct sum of endomorphism algebras for highest weight vectors. We then describe in general terms the structure of a Lie superalgebra preserving an odd Hermitean nondegenerate inner product, like the Poincar\'e one that we use. In Section \ref{sec:evenoperators} we first build an inclusion of $\mathbf{sl}(4,\R)$ inside $\ggr$, and then use three toral operators $K_j$ built as supercommutators of $V_j$ and $A_j$ to decompose the generators of this 
$\mathbf{sl}(4,\R)$ into $K_j$-homogeneous components. In Section \ref{sec:bases} we describe explicitely some naturally structured bases for the spaces of highest weight vectors inside the isotypical components mentioned above, and in Section \ref{sec:restriction} we restrict $\ggr$ to these spaces of highest weight vectors , and use the bases previously built and the information accumulated up to this point to obtain the final description of the structure of $\ggr$ and of $\gg$.
In section \ref{sec:finalremarks} we show that both $\ggr$ and $\gg$ are closed with respect to the "standard" adjunction operator induced by the metric, and therefore $\gg$ is a $*$-Lie superalgebra.\\
The tables in the Appendix describe in detail the action of the $K_j$-homogeneous components of the generators of the $\mathbf{sl}(4,\C)$ mentioned above, in the explicitely provided bases.

\section{The algebra and its invariance properties}
\label{sec:alginv}
Let us fix a point $p$ in the
WSD manifold $X$. In accordance with the similar notation introduced in \cite{GG} for the rank two case, the WSD structure splits the cotangent space as
$T_p^*X = W_0 \oplus W_1 \oplus W_2$ where the $W_j$ are  three
mutually orthogonal three-dimensional canonical distributions defined as:
\[W_0 = \{\phi\in T_p^*X~|~\phi\wedge \omega_1^3 = \phi\wedge
\omega_2^3 = 0\}\]
\[W_1 = \{\phi\in T_p^*X~|~\phi\wedge \omega_1^3 = \phi\wedge
\omega_D^3 = 0\}\]
\[W_2 = \{\phi\in T_p^*X~|~\phi\wedge \omega_2^3 = \phi\wedge
\omega_D^3 = 0\}\]

\begin{dfn}
The decomposition $T_p^*X = W_0 \oplus W_1 \oplus W_2$ induces corresponding decompositions for $s\in\{0,...,9\}$:
\[\bigwedge ^s T_p^*X = \bigoplus_{p+q+r = s}\bigwedge^p W_0 \oplus \bigwedge^q W_1 \oplus \bigwedge^r W_2\]
We say that $\alpha\in \bigwedge^p W_0 \oplus \bigwedge^q W_1 \oplus \bigwedge^r W_2$ has {\em multidegree} $(p,q,r)$
\end{dfn} 
The WSD structure also determines canonical pairwise linear
identifications among $W_0,W_1$ and $W_2$, so that one can also
write $T_p^*X = W_0\otimes_\R \R^3$ or more simply
\[T_p^*X = W\otimes_\R \R^\dime\] where $W = W_0 \cong W_1 \cong
W_2$.\\
The canonical decomposition above equips $W$  (and hence on all the $W_j$) with an orientation  induced by the one on $T^*_pX$.
Let us now come back to the canonical  operators $L_j,V_j$ mentioned in
the Introduction using the notation $Vol(W_j)$ for the oriented volume forms of the $W_j$: \\
\textbf{Definition~\ref{dfn:Loperators}}\textit{ For $\phi\in
\Omega^*_\C X$,
\[L_0(\phi) = \omega_D\wedge \phi,\qquad L_1(\phi) = -\omega_2\wedge
\phi,\qquad L_2(\phi) = \omega_1\wedge \phi\] 
\[V_0(\phi) = Vol(W_0)\wedge \phi,\qquad V_1(\phi) = Vol(W_1)\wedge \phi, \qquad V_2(\phi) = Vol(W_2)\wedge \phi\] }\\
We now choose a (non-canonical) orthonormal basis
$\gamma_1,\gamma_2,\gamma_3$ for $W_0$, and this together with the standard
identifications of the $W_j$ determines an orthonormal basis for
$T_p^*X$, which we write as $\{v_{ij}= \gamma_i\otimes
e_j~|~i=1,2,3,~j=0,1,2\}$.  We remark that the $v_{ij}$ are an
\textit{adapted coframe} for the WSD structure, and therefore we
have the explicit expressions:
\[\omega_1 = v_{10}\wedge v_{11} + v_{20}\wedge v_{21} + v_{30}\wedge v_{31}\]
\[\omega_2 = v_{10}\wedge v_{12} + v_{20}\wedge v_{22} + v_{30}\wedge v_{32}\]
\[\omega_D = v_{11}\wedge v_{12} + v_{21}\wedge v_{22} + v_{31}\wedge v_{32}\]
A different choice of the $\gamma_1,\gamma_2,\gamma_3$ would be related to
the previous one by an element in $\mathbf{O}(3,\R)$ or, taking into
account the orientability of $X$ mentioned in the Introduction, an
element of $\mathbf{SO}(3,\R)$. The Lie algebra of the group
$\mathbf{SO}(3,\R)$ expressing the change from one oriented adapted
basis to another is generated (point by point) by the global
operators $J_1,J_2,J_3$:
\begin{dfn}
\label{dfn:so2R} For fixed $p\in X$, the operators $J_1,J_2,J_3 \in End_\R(\Lambda^*T^*_pX)$ are defined in terms of the standard basis $v_{ij}$ as
\[J_1(v_{2j}) = v_{3j},\quad J_1(v_{3j}) = - v_{2j},\quad J_1(v_{1j}) = 0\qquad \text{for}~
j \in \{0,1,2\}\]
\[J_2(v_{3j}) = v_{1j},\quad J_2(v_{1j}) = - v_{3j},\quad J_2(v_{2j}) = 0\qquad \text{for}~
j \in \{0,1,2\}\]
\[J_3(v_{1j}) = v_{2j},\quad J_3(v_{2j}) = - v_{1j},\quad J_3(v_{3j}) = 0\qquad \text{for}~
j \in \{0,1,2\}\]
and $J_i(v\wedge w) = J_i(v)\wedge w + v \wedge J_i(w)$
for $v,w \in \Lambda^*T^*_pX$ and $i\in\{1,2,3\}$.
\end{dfn}
\begin{rmk} The operators $J_1,J_2,J_3$ generate  the Lie algebra $\mathbf{aut}_p(X)$ of all the operators preserving the WSD structure at the given point, which does not depend on the choice of a basis for $T^*_pX$. If the metric of $X$ is such that the structural forms $\omega_1,\omega_2,\omega_D$ are covariant constant, then its holonomy must lie inside the group  $\mathbf{Aut}_p(X) = exp(\mathbf{aut}_p(X))$.
\end{rmk}

Using the chosen (non canonical) orthonormal basis, one can define corresponding wedge and contraction operators:
\begin{dfn} Let $i\in\{1,2,3\}$ and $j\in\{0,1,2\}$. The operators $E_{ij}$ and $I_{ij}$ are respectively the wedge and the contraction operator with the form $v_{ij}$ on $\bigwedge^*T^*X$ (defined using the given basis); we use the notation $\frac{\partial}{\partial v_{ij}}$ to indicate the element of $T_pX$ dual to $v_{ij}\in T^*_pX$:
\[E_{ij}(\phi) = v_{ij}\wedge \phi,\qquad I_{ij}(\phi) = \frac{\partial}{\partial v_{ij}}\contr
\phi\]
\end{dfn}
\begin{pro}
\label{pro:clifrelations}
The operators $E_{ij},I_{ij}$ satisfy the following relations:
\[ \forall i,j,k,l\qquad E_{ij}E_{kl} = -E_{kl}E_{ij},\quad I_{ij}I_{kl} = - I_{kl}I_{ij}\]
\[\forall i,j \qquad E_{ij}I_{ij} + I_{ij}E_{ij} = Id\]
\[\forall (i,j)\not= (k,l) \qquad E_{ij}I_{kl} = - I_{kl}E_{ij}\]
\[\forall i,j\qquad E_{ij}^* = I_{ij},\quad I_{ij}^* = E_{ij}\]
where $*$ is adjunction with respect to the metric.
\end{pro}
\dimo The proof is a simple direct verification, which we omit. \qed\\

Using the (non canonical) operators $E_{ij}$ we can obtain simple expressions for the
pointwise action of the other canonical operators, the odd operators $V_j$ induced by the volume forms of the distributions $W_j$:
 for $\phi\in \bigwedge^*T^*_pX$,
\[V_0(\phi) = E_{10} E_{20}E_{30}(\phi),\qquad
V_1(\phi) = E_{11} E_{21}E_{31}(\phi),\qquad V_2(\phi) =
E_{12} E_{22}E_{32}(\phi)\]
Remember however that the operators $V_j$ do not depend on the choice of a basis, as they
are simply multiplication by the volume forms of the spaces $W_j$.\\
We use the $v_{ij}$ also as a orthonormal basis for the complexified
space $T_p^*\otimes_\R\C$ (with respect to the induced hermitian
inner product). We indicate with the same symbols $V_j$ the
complexified operators acting on the spaces $\bigwedge_\C^*T^*_pX$. The riemannian metric induces a riemannian metric on $T^*_p X$ and
on the space $\bigwedge^*T^*_pX$. 
\begin{rmk}
The Riemannian metric on $X$ induces in the standard way a Hodge star operator
\[*: {\bigwedge}_\C^*T^*_pX\to {\bigwedge}_\C^*T^*_pX\]
which (as $dim(X) = 9$) satisfies $*^2 = Id$. 
\end{rmk}
As usual, one can  define the induced Hermitean inner product $(~,~)_p$ on $\bigwedge_\C^*T^*_pX$ via  the Hodge $*$ operator
\begin{dfn}
Given a homogeneous operator $\phi$ on $\bigwedge_\C^*T^*_pX$, we indicate with $\phi^\ss$ the homogeneous operator which satisfies:
\[ (\phi(\alpha),\beta)_p = (-1)^{deg(\phi)deg(\alpha)}(\alpha,\phi^\ss(\beta))_p\]
for all homogeneous $\alpha,\beta$.
\end{dfn}
\begin{dfn}
\label{dfn:lambda-A-operators} For $j\in\{0,1,2\}$
\[\Lambda_j = L_j^\ss,\qquad A_j =
V_j^\ss\]
\end{dfn}

By construction the canonical operators $L_j,V_j,\Lambda_j,A_j$ on
$\bigwedge_\C^*T^*_pX$ are the pointwise restrictions of
corresponding global operators on smooth differential forms, which
we indicate with the same symbols: for $j \in \{0,1,2\}$,
\[L_j,V_j,\Lambda_j,A_j: \Omega_\C^*(X)\to
\Omega_C^*(X)\]
Summing up:
\begin{dfn}
\label{dfn:algebragg} The Lie superalgebra $\ggr$ is the Lie
subalgebra of the general linear Lie superalgebra of $\Omega^*(X)$ generated by the
operators
\[\{\imath L_j,\imath V_j,\imath \Lambda_j,A_j~|~\text{for}~ ~ j = 0,1,2\}\]
The Lie superalgebra $\gg$ is $\ggr\otimes \C$, and is in a natural way a Lie subalgebra of the general linear Lie superalgebra of $\Omega_\C^*(X)$.
\end{dfn}
\begin{dfn}
\label{dfn:superhermitean}
For every $p\in X$ there is a natural odd non degenerate super Hermitean inner product $<~,~>_p$ on $\bigwedge_\C^*T^*_pX$, defined using the natural (standard) Hermitean inner product $(~,~)_p$ and the (pointwise) volume form $\Omega$ associated to the Riemannian metric and the orientation:
\[ <\alpha,\beta>_p = (\alpha\wedge \overline{\beta},\Omega)_p\]
\end{dfn}
We list without proof the following standard facts:
\begin{pro} For every $p\in X$, the pairing $<~,~>_p$ satisfies the following properties:\\

a) $\quad < \alpha,\beta>_p = (\alpha, *\beta)_p$\\

b) $\quad <\alpha,\beta> = (-1)^{deg(\alpha)deg(\beta)}\overline{<\beta,\alpha>_p}\quad\text{(super Hermitean)}$\\

c) $\quad \forall\alpha\left(\forall\beta <\alpha,\beta>_p = 0\implies \alpha = 0\right)\quad\text{(non-degenerate)}$\\

d) $\quad <~,~>_p$ is  preserved by the action of the group of orientation preserving isometries $SO(T^*X_p,\me)$ on $\bigwedge_\C^*T^*_pX$ (and in particular by our $\mathbf{SO}(3,\R)$).\\

e) $\quad <~,~>_p$ is  preserved by the Hodge $*$ operator.
\end{pro}
\begin{rmk}
As $dim(X)$ is odd, for $\alpha,\beta$ homogeneous $<\alpha,\beta>\not=0$ implies that $\{deg(\alpha),\deg(\beta)\}$ contains an even and an odd number, and therefore the product $<~,~>_p$ is actually \textit{Hermitean}, besides being super Hermitean.
\end{rmk}
The choice of the particular real form $\ggr$ is motivated by the following theorem:
\begin{teo}
\label{teo:superantihermitean} The algebra $\ggr$  preserves the  form $<~,~>_p$ for every $p\in X$. This is equivalent to the fact that the condition $\phi^\ss = -*\,\phi\,*$ must hold for all the elements $\phi$ of $\ggr$.
\end{teo}
\dimo
For fixed $p\in X$, by definition, a  homogeneous operator $\phi$ preserves $<~,~>_p$ if and only if it satisfies the following equation for homogeneous $\alpha,\beta$:
\[\forall\alpha\forall\beta \quad <\phi(\alpha),\beta>_p + (-1)^{deg(\alpha)deg(\phi)}<\alpha,\phi(\beta)>_p = 0\]
Using the expression $<\alpha,\beta>_p  = (\alpha,*\beta)_p$, the above is equivalent to $\phi^\ss = -*\,\phi\,*$.\\
It is enough to check the equation for the generators of $\ggr$, as one has that if $\phi$ and $\psi$ are homogeneous and satisfy the equation, then also :
\[<[\phi,\psi](\alpha),\beta>_p + (-1)^{deg(\alpha)(deg(\phi)+deg(\psi))}<\alpha,[\phi,\psi](\beta)>_p = 0\]
(by general reasoning in the category of super vector spaces and operators, or by direct computation).
For the generators $\imath L_j$ and $\imath V_j$, which are wedge operators with homogeneous differential forms, the verification of the equation is immediate. The operators $i\Lambda_j$ by definition satisfy the equation
\[\forall \alpha,\beta\quad (\imath L_j(\alpha),\beta)_p = -(\alpha,\imath \Lambda_j (\beta))_p\]
Using the Hodge $*$, one can rephrase the previos equation as:
\[\forall \alpha,\beta\quad\imath L_j(\alpha)\wedge \overline{*\beta} = -\alpha\wedge \overline{*\imath\Lambda_j(\beta)}\]
As $L_j$ is a wedge operator with an even form, from this one obtains
\[\imath \Lambda_j = * ~ \imath L_j ~ *\]

As the Hodge $*$ is self-adjoint with respect to the product $<~,~>_p$ and $\imath L_j$ is anti-(super)selfadjoint, one has
\[<\imath\Lambda_j(\alpha),\beta>_p = <*~\imath L_j~*(\alpha),\beta>_p = - < \alpha, *~\imath L_j~*(\beta)>_p = -<\alpha,\imath\Lambda_j(\beta)>_p\] 
It remains to verify the case of $A_j$. 
 The operators $A_j$ by definition satisfy the equation
\[\forall \alpha,\beta\quad (V_j(\alpha),\beta)_p = (-1)^{deg(\alpha)}(\alpha,A_j (\beta))_p\]
Using the Hodge $*$, one can rephrase the previos equation as:
\[\forall \alpha,\beta\quad V_j(\alpha)\wedge \overline{*\beta} = (-1)^{deg(\alpha)}\alpha\wedge \overline{*A_j(\beta)}\]
As $V_j$ is a wedge operator with an odd form, from this one obtains
\[A_j = * ~ V_j ~ *\]
As before, we start with the observation that $V_j$ is superselfadjoint and the Hodge star is selfadjoint (with respect to $<~,~>_p$). Therefore
\[<A_j(\alpha),\beta>_p = <*~ V_j~*(\alpha),\beta>_p = \]
\[= (-1)^{deg(*\alpha)}< \alpha, *~V_j~*(\beta)>_p = -(-1)^{deg(\alpha)}<\alpha,A_j(\beta)>_p\] 
\qed
\section{The actions of  $\mathcal{S}_3$ and $\so3$}
\label{sec:so2R}

The canonical splitting
$T_p^*X = W_0 \oplus W_1 \oplus W_2$ together with the canonical
identifications $W_0\cong W_1\cong W_2$ induce an action of the
symmetric group $\mathcal{S}_3$ on $T^*_pX$, which propagates to
$\bigwedge^*T^*X$ and to its $\mathcal{C}^\infty$ sections. At
every point, the action can be written explicitly in terms of the
basis $\{v_{ij}\}$ from the beginning of section \ref{sec:alginv} as
\[\sigma(v_{ij}) = v_{i\sigma(j)}\]
The induced action on endomorphisms via
conjugation, $\sigma(\phi) = \sigma\circ\phi\circ\sigma^{-1}$,
preserves $\ggr$. Indeed, one can check directly using the basis
$v_{ij}$ at every point that for $\sigma \in \s3$
\[\sigma(V_j) = V_{\sigma(j)},\qquad \sigma(L_j) = \epsilon(\sigma)L_{\sigma(j)}\]
Since $\s3$ acts on $\ggr$ by conjugation with (even) unitary operators, its
action commutes with adjunction (the $\ss$ operator), and therefore
\[\sigma(A_j) = A_{\sigma(j)},\qquad \sigma(\Lambda_j) = \epsilon(\sigma)\Lambda_{\sigma(j)}\]
Moreover, one immediately verifies that $\sigma(J_i) = J_i$  for all $\sigma\in\s3$ and for all $i\in\{1,2,3\}$, which means that the
action of $\mathcal{S}_3$ commutes with that of $\so3$.

\begin{teo}
\label{invso3}
The algebra $\ggr$ commutes with the action of $\mathbf{SO}(3,\R)$.
\end{teo}
\dimo
It is enough to show that the generators  of $\ggr$ commute with the  (even)  generators $J_1,J_2,J_3$ of $\mathbf{so}(3,\R)$.
As the operators $J_k$ are (real) antiHermitean with respect to the standard scalar product $(~,~)_p$, it is enought to check the commutation of them with the $L_j$ and the $V_j$ (indicating with $Vol\,(W_j)$ the oriented volume form of the distribution $W_j$).
\[[L_j,J_k](\alpha) = \omega_j\wedge J_k(\alpha) - J_k(\omega_j\wedge \alpha) =  - J_k(\omega_j) \wedge \alpha\]
\[[\imath V_j,J_k](\alpha) = \imath Vol\,(W_j)\wedge J_k(\alpha) - J_k(\imath Vol\,(W_j))\wedge \alpha) =  - J_k(\imath Vol\,(W_j)) \wedge \alpha\]
As the $J_k$ are antiHermitean and preserve the distributions $W_j$,  $J_k(Vol\,(W_j)) = 0$. To verify that $J_k(\omega_j) = 0$ one makes a direct computation in (orthonormal) coordinates using Definition \ref{dfn:so2R}. Using the action of $\mathcal{S}_3$, one actually reduces to the following two verifications:
\[J_1(\omega_1) = 0,\qquad J_1(\omega_2) = 0\]
\[J_1(\omega_1) = J_1(v_{10}\wedge v_{11} + v_{20}\wedge v_{21} + v_{30}\wedge v_{31}) =\]
\[= v_{30}\wedge v_{21} + v_{20}\wedge v_{31} - v_{20}\wedge v_{31} - v_{30}\wedge v_{21} = 0\] 
\[J_1(\omega_2) = J_1(v_{10}\wedge v_{12} + v_{20}\wedge v_{22} + v_{30}\wedge v_{32}) =\]
\[= -v_{30}\wedge v_{12} - v_{10}\wedge v_{32} + v_{10}\wedge v_{32} + v_{30}\wedge v_{12} = 0\]
\qed
\begin{rmk} The algebra $\ggr$ (and its complexification $\gg$) commutes with the complexification $\mathbf{so}(3,\C)$ of the Lie algebra generated by the operators $J_k$.
\end{rmk}
\begin{dfn} Given $n\in\Z_{\geq 0}$, we indicate with $\rho_n$ the $2n +1$
dimensional complex irreducible representation of  the Lie algebra $\mathbf{so}(3,\C)$, which is a highest weight representation with highest weight equal to $2n$.
\end{dfn}
The actual representation of $\mathbf{so}(3,\R)$ (and $\mathbf{so}(3,\C)$) on $\bigwedge^*T^*_\C(X)$ does depend on the choice of an orthonormal basis (indeed, two different adapted oriented orthonormal bases differ by an element in $\mathbf{SO}(3,\R)$). However, the decomposition in isotypical components described in the following proposition is canonical, and does not depend on any choice.
\begin{pro}
\label{pro:table} Under the $\mathbf{so}(3,\C)$ representation induced by the
operators
$J_k$, for any $p\in X$ :\\
a) The space $T^*_\C(X)_p$ decomposes as an  $\mathbf{so}(3,\C)$ module as
\[T^*_\C(X)_p = \left(W_0\otimes\C\right) \oplus \left(W_1\otimes\C\right) \oplus \left(W_2\otimes\C\right) \]
with $W_j\otimes\C \cong  \rho_1$ as  $\mathbf{so}(3,\C)$ modules.\\
b) More generally, Table \ref{tab:isotyp} describes the decomposition of $\bigwedge^*T^*_\C(X)_p$ into isotypical components
\begin{table}[hbtp]
\begin{center}
\begin{tabular}{|l|c|c|c|c|}
\hline
  & &  &       & \\
    & Type $\rho_0$ & Type $\rho_1$ & Type  $\rho_2$ & Type $\rho_3$ \\
      & &  &       & \\
\hline
  & &  &       & \\
$\bigwedge^0T^*_\C X$ &  $\rho_0$   &  &       & \\
  & &  &       & \\
$\bigwedge^1T^*_\C X$ &     &  $\rho_1^{\oplus 3}$ &       & \\
  & &  &       & \\
  $\bigwedge^2T^*_\C X$ & \(\rho_0^{\oplus 3}\)    &  $\rho_1^{\oplus 6}$ &    \(\rho_2^{\oplus 3}\)    &  \\
  & &  &       & \\
  $\bigwedge^3T^*_\C X$ &   \(\rho_0^{\oplus 10} \) &  $\rho_1^{\oplus 9}$ &  \(  \rho_2^{\oplus 8} \)   &  \(\rho_3\)\\
  & &  &       & \\
  $\bigwedge^4T^*_\C X$ &  \( \rho_0^{\oplus 6} \) &  $\rho_1^{\oplus 18}$ & \(  \rho_2^{\oplus 9}   \)  & \(\rho_3^{\oplus 3}\) \\
  & &  &       & \\
  $\bigwedge^5T^*_\C X$ &\(\rho_0^{\oplus 6}  \)  &  $\rho_1^{\oplus 18}$ &\(  \rho_2^{\oplus 9}   \)   & \(\rho_3^{\oplus 3}\) \\
  & &  &       & \\
  $\bigwedge^6T^*_\C X$ &\(  \rho_0^{\oplus 10} \)  &  $\rho_1^{\oplus 9}$ &\(    \rho_2^{\oplus 8}   \) & \(\rho_3\) \\
  & &  &       & \\
  $\bigwedge^7T^*_\C X$ &\( \rho_0^{\oplus 3} \)   &  $\rho_1^{\oplus 6}$ &     \(\rho_2^{\oplus 3}\)   &  \\
  & &  &       & \\
  $\bigwedge^8T^*_\C X$ &     &  $\rho_1^{\oplus 3}$ &       & \\
  & &  &       & \\
  $\bigwedge^9 T^*_\C X$ & \( \rho_0\)   &   &       & \\
    & &  &       & \\
\hline
 & &  &       & \\
  Dimension & $40$   & $3\times 72$  & $5\times 40$     & $7\times 8$\\
    & &  &       & \\
\hline
\end{tabular}
\medskip
\caption{Isotypical components} \label{tab:isotyp}
\end{center}
\end{table}
\end{pro}
\dimo a) The subspaces $W_j$ are preserved by the action of $\mathbf{so}(3,\R)$, therefore there is a direct sum decomposition of $T^*_\C(X)_p$ with respect to the action of $\mathbf{so}(3,\C)$. \\
The action of $\mathbf{so}(3,\C)\cong \mathbf{sl}(2,\C)$ on $\bigwedge^*T^*_\C X$ is best explained using a Serre set of generators $\mathbf{e},\mathbf{f},\mathbf{h}$, defined as follows in terms of the operators of Definition \ref{dfn:so2R}:
\[\mathbf{e} = \imath J_1 - J_2,\quad \mathbf{f} = \imath J_1 + J_2,\quad \mathbf{h} = 2\imath J_3\]
To show the isomorphisms with $\rho_1$ we use a new orthonormal basis made of eigenvectors for $\mathbf{h}$, built starting from the basis $\{v_{ij}\}$:
\[w_{1j} = v_{1j} + \imath \,v_{2j},\quad w_{2j} = v_{1j} - \imath \,v_{2j},\quad w_{3j} = v_{3j}\]
One then checks that for fixed $j$ the space $<w_{1j},w_{2j},w_{3j}>_\C = W_j\otimes \C$ isomorphic as $\mathbf{sl}(2,\C)$-module to $\rho_1$, since  $w_{1j}$ is a highest weight vector of $\mathbf{h}$-weight $2$.\\
b) The proof consists now of some plethysm for $\mathbf{sl}(2,\C)$. We will use the well known formula
\[\rho_a\otimes \rho_b = \sum_{i=|a-b|}^{a+b}\rho_i\]
and the facts that $\bigwedge^2\rho_1 \cong \rho_1$ and $\bigwedge^3\rho_1 \cong \rho_0$.\\
The space $\bigwedge^2T^*_\C X$ is isomorphic as an $\mathbf{sl}(2,\C)$ module to $\bigwedge^2\left(\rho_1^\ot\right)$, and therefore
\[\bigwedge^2T^*_\C X \cong \bigwedge^2\left(\rho_1^\ot\right) \cong \left(\bigwedge^2\rho_1\right)^\ot + \left(\rho_1\otimes\rho_1\right)^\ot \cong \rho_1^\ot + \left(\rho_0 + \rho_1 + \rho_2\right)^\ot\]
The space $\bigwedge^3T^*_\C X$ is isomorphic as an $\mathbf{sl}(2,\C)$ module to $\bigwedge^3\left(\rho_1^\ot\right)$, and therefore
\[\bigwedge^3T^*_\C X \cong 
\bigwedge^3\left(\rho_1^\ot\right)
\cong 
\left(\bigwedge^3\rho_1\right)^\ot + \left[\left(\bigwedge^2\rho_1\right)\otimes 
\rho_1\right]^\os + 
\quad \rho_1\otimes\rho_1\otimes\rho_1 \cong \]
\[\cong \rho_0^\ot + \left(\rho_2 + \rho_1 + \rho_0\right)^\os + \left(\rho_2 + \rho_1 + \rho_0\right)\otimes\rho_1 \cong \rho_0^{\oplus 10} + \rho_1^{\oplus 9} + \rho_2^{\oplus 8} + \rho_3 \]
The space $\bigwedge^4T^*_\C X$ is isomorphic as an $\mathbf{sl}(2,\C)$ module to $\bigwedge^4\left(\rho_1^\ot\right)$, and therefore
\[\bigwedge^4T^*_\C X \cong \bigwedge^4\left(\rho_1^\ot\right) \cong
\left[\left(\bigwedge^3\rho_1\right)\otimes \rho_1\right]^\os +
\left[\left(\bigwedge^2\rho_1\right)\otimes \left(\bigwedge^2\rho_1\right)\right]^\ot +\]
\[+  \left[\left(\bigwedge^2\rho_1\right)\otimes\rho_1\otimes\rho_1\right]^\ot\cong \left[\rho_1\right]^\os + \left[\rho_0 + \rho_1 + \rho_2\right]^\ot +  \left[\rho_0 + \rho_1^{\oplus 3} + \rho_2^{\oplus 2} + \rho_3\right]^\ot
 \]
As the Hodge $*$ commutes with the induced action of the isometries of $T_pX$, it must preserve the isotypical components with respect to $\mathbf{so}(3,\C)$. This implies that the computations above are enough to determine the remaining decompositions, as one has the  isomorphisms of  $\mathbf{so}(3,\C)$ modules $*:\bigwedge^kT^*_\C X \to \bigwedge^{9-k}T^*_\C X$.
\qed\\

Let us call $HW_k$ the complex vector space of $\mathbf{so}(3,\C)$-highest weight vectors inside the isotypical component of type $\rho_k$.
As the algebra $\ggr$ commutes with $\mathbf{so}(3,\C)$, it preserves each isotypical component, and moreover it sends each space $HW_k$ to itself.
\begin{cor}
\label{cor:subend}
The algebra $\ggr$ is mapped injectively to the space
\[End_\C(HW_0)\oplus End_\C(HW_1)\oplus End_\C(HW_2)\oplus End_\C(HW_3)\]
and therefore its dimension over $\R$ is at most $2\times (40^2 + 72^2 + 40^2 + 8^2) = 16896$
\end{cor}
\dimo It only remains to check the dimensions of the $HW_k$, which can be read from  Table \ref{tab:isotyp}. \qed\\
As a consequence of this corollary, the upper bound on the dimension over $\R$ of the algebra $\ggr$ has been reduced from $2^{19} = 524288$ to $16896$. To further reduce it, we make the following observation:
\begin{pro}
\label{pro:injectively}
The odd nondegenerate super Hermitean form $<~,~>_p$ remains nondegenerate when restricted to the spaces $HW_k$. The algebra $\ggr$ is mapped injectively to the space of operators which are super antiHermitean  for this form and have supertrace zero:
\[ \mathbf{su}(HW_0^{ev}|HW_0^{odd},<~,~>_p)~ \oplus~  \mathbf{su}(HW_1^{ev}|HW_1^{odd},<~,~>_p)~ \oplus\]
\[ \oplus ~ \mathbf{su}(HW_2^{ev}|HW_2^{odd},<~,~>_p)~ \oplus~ \mathbf{su}(HW_3^{ev}|HW_3^{odd},<~,~>_p)\]
\end{pro}
\dimo Let $\alpha$ be a nonzero element of $HW_k$. As the Hodge $*$ commutes with the action of $\mathbf{so}(3,\C)$, $*\alpha$ is again a (nonzero) element of $HW_k$. Moreover,
\[<\alpha,*\alpha>_p = (\alpha,\alpha)_p > 0\] 
The second claim is a consequence of Theorem ~\ref{teo:superantihermitean}. \\
The supertrace is zero because all the generators have supertrace equal to zero (the even generators $\imath L_j,\imath \Lambda_j$ are nilpotent). 
\qed\\

As the form $<~,~>$ is odd and nondegenerate, the even and odd parts of the space must be in duality (like in the case of the $HW_k$). If we fix a  Hermitean inner product on the odd part (and an induced compatible one on the even part), this duality can be turned into a linear isometry which we indicate with $*$, as in our case it is the Hodge $*$. As the manifold $X$ has odd dimension, the Hodge $*$ has also the property that $*^2 = Id$, so that this will hold also in all the $HW_k$. Summing up, there is an orthogonal decomposition
\[HW_k ~=~ HW_k^{ev} \oplus^\perp HW_k^{odd}\]
and the odd pairing $<~,~>_p$ puts in duality the odd part with the even one. The Hodge $*$ gives an isometric involution from $HW_k$ to itself wich exchanges $HW_k^{ev}$ with   $HW_k^{odd}$.\\
The standard arguments above can be summarized in the  explicit description of all the maps which preserve the odd form $<~,~>_p$:
\begin{pro}
\label{pro:sync}
Let $\phi: HW_k\to HW_k$ be a linear morphism whose homogeneous components preserve the odd pairing $<~,~>_p$ (i.e. $\phi$ preserves $<~,~>_p$). With respect to the above mentioned  decomposition and isomorphism 
\[HW_k ~=~ HW_k^{ev} \oplus^\perp HW_k^{odd},\quad HW_k^{ev}\cong \left(HW_k^{odd}\right)^*\]
$\phi$ must satisfy $\phi^\ss = -*\,\phi\,*$ and therefore it has the following form (indicating with $\phi_j^*$ the adjoint with respect to the Hermitean inner product)
\[\phi = \begin{pmatrix}\phi_1 & \phi_2 = \phi_2^*\\ \phi_3 = -\phi_3^* & -\phi_1^*\end{pmatrix}\]
where $\phi_1:HW_k^{ev}\to HW_k^{ev}$, $\phi_2:HW_k^{odd}\to HW_k^{ev}$, $\phi_3:HW_k^{ev}\to HW_k^{odd}$.
\end{pro}
\dimo The condition follows immediately from $\phi^\ss = -*\,\phi\,*$ applied to the "matrix" of components of the operator $\phi$. \qed\\
From the previous propositions and the fact that the supertrace of the elements of $\ggr$ is zero one has
\begin{cor} The algebra $\ggr$ has dimension at most 
\[ 4*20^2 + 4*36^2 + 4*20^2 + 4*4^2 -4 =  8444\]
\end{cor}
In the remainder of the paper we will show that this estimate is almost sharp, being off by $48$, and the previous propositions will be turned into a precise description of the algebra $\ggr$ (see Theorem \ref{teo:mainteo}). 
\section{Some even operators and global relations}
\label{sec:evenoperators}
In this section we study some even operators obtained as superbrackets of the odd generators, the $K_{lm}$, and study how they relate to the even algebra generated by the $\imath L_j,\imath\Lambda_j$. 
\begin{dfn} For $j \in \{0,1,2\}$
\[H_j = [\imath\Lambda_j,\imath L_j]\]
\end{dfn} 
\begin{pro}
The operators $\{\imath L_0,\imath \Lambda_0, \imath L_1,\imath \Lambda_1, \imath L_2,\imath \Lambda_2\}$ generate the real Lie algebra $\mathbf{sl}(4,\R)$, and the  $\{H_0,H_1,H_2\}$ span a Cartan subalgebra.
\end{pro}
\dimo Referring to \cite{GG},\cite{G2}, one sees that there is a Serre basis for the complexified algebra, given by the operators
\[\begin{array}{lll}\mathbf{e}_1 = [L_0,\Lambda_1] &  \mathbf{f}_1 = [L_1,\Lambda_0] & \mathbf{h}_1 = [\mathbf{e}_1,\mathbf{f}_1]\\
\mathbf{e}_2 = [L_1,\Lambda_2] &  \mathbf{f}_2 = [L_2,\Lambda_1] & \mathbf{h}_2 = [\mathbf{e}_2,\mathbf{f}_2]\\
\mathbf{e}_3 = L_2 &  \mathbf{f}_3 = \Lambda_2 & \mathbf{h}_3 = [\mathbf{e}_3,\mathbf{f}_3]\end{array}\]
which is associated to the Dynkin diagram $A_3$. From this it is immediate to check that the real Lie algebra generated by $H_j,\imath L_j,\imath \Lambda_j$ is a split form,  isomorphic to $\mathbf{sl}(4,\R)$. This proof carries over also in arbitrary rank (of the WSD structure).\qed\\

In the rank $2$ case (see \cite{GG}) one verified that the brakets of the form $[V_j,A_m]$ were already contained in the Cartan subalgebra spanned by the $H_k$. Contrary to that situation, in rank $3$ the algebra is richer, and we get genuinely new operators:
\begin{dfn} For $l,m \in \{0,1,2\}$
\[K_{lm} = \imath V_lA_m + \imath A_m V_l = [\imath V_l,A_m]\]
\end{dfn} 
\begin{pro} For $j,l,m \in \{0,1,2\}$ one has
\[[H_j,K_{lm}] = \left(-3\delta_{jl} + 3\delta_{jm}\right)\, K_{lm}\]
In particular, $[H_j,K_{mm}] = 0$ for all $j,m \in \{0,1,2\}$.
\end{pro}
\dimo This follows using (super)Jacoby from the relations
\[[H_j,V_m] = 3(1-\delta_{jm})V_m,\quad [H_j,A_m] = -3(1-\delta_{jm})A_m\]
which can be computed using the expressions in terms of the $E_{ij}$, as in \cite{GG}.\qed\\
\begin{dfn}
We say that $\alpha\in \bigwedge^*T^*_p X$ has $Kw$ weight equal to $(z_0,z_1,z_2)$ if 
\[K_{mm}(\alpha) = z_m\alpha\]
for every $m\in\{0,1,2\}$
\end{dfn}

The $H_j,K_{mm}$ form an abelian subalgebra from the previous proposition. However, although the generators of $\ggr$ are already eigenvectors for the adjoint action of the $H_j$, they are not so for the adjoint action of the $K_{mm}$. Nevertheless, the adjoint action of the $K_{mm}$ can be diagonalized (over the real numbers) on the Lie subalgebra of $\ggr$ generated by the $\imath L_j$ and the $\imath \Lambda_k$:
\begin{pro}
\label{pro:kdegree}
If a form $\alpha\in \bigwedge^*T^*_p X$ is multidegree homogeneous of multidegree $(a,b,c)$ then is also $Kw$ homogeneous, with $Kw$-weight:
\[Kw(\alpha) = \imath(-1)^{deg(\alpha)}(\delta_{a0} - \delta_{a3},\delta_{b0} - \delta_{b3},\delta_{c0} - \delta_{c3})\]
As a consequence, the operators $K_{mm}$ can be silmultaneously diagonalized over $\C$, with eigenvalues $0$ and $\pm \imath$.
\end{pro}
\dimo
It is enough to prove the proposition for $\alpha = v_{i_1 j_1}\wedge\cdots \wedge v_{i_k j_k}$ a monomial in the forms $v_{ij}$, which consitute as a basis for $\bigwedge^*T^*_pX$. Such a monomial of multidegree $(a,b,c)$  can be written as $m_0\wedge m_1 \wedge m_2$, where $m_0$ has multidegree $(a,0,0)$, $m_1$ has multidegree $(0,b,0)$ and $m_2$ has multidegree $(0,0,c)$. The action of $\imath V_m$ and of $A_m$ on such a monomial is easily described:
\[\imath V_0\left(m_0\wedge m_1 \wedge m_2\right) = \imath \delta_{a0}Vol(W_0)\wedge m_1\wedge m_2\]
\[A_0\left(m_0\wedge m_1 \wedge m_2\right) = \delta_{a3}(-1)^{b+c} m_1 \wedge m_2\]
\[\imath V_1\left(m_0\wedge m_1 \wedge m_2\right) = \imath \delta_{b0}Vol(W_1)\wedge m_0\wedge m_2\]
\[A_1\left(m_0\wedge m_1 \wedge m_2\right) = \delta_{b3}(-1)^{c} m_0 \wedge m_2\]
\[\imath V_2\left(m_0\wedge m_1 \wedge m_2\right) = \imath \delta_{c0}Vol(W_2)\wedge m_0\wedge m_1\]
\[A_2\left(m_0\wedge m_1 \wedge m_2\right) = \delta_{c3}\, m_0 \wedge m_1\]
From these explicit expressiont it is straightforward to verify that in this basis the $K_{mm} = \imath V_m A_m + \imath A_m V_m$ are diagonal and their eigenvalues   are $0$ and $\pm \imath$. 
More precisely, on a monomial $m$ of multidegre $(a,b,c)$ one has therefore that
\[K_{00}(m) = \imath(-1)^{deg(m)}(\delta_{a0} - \delta_{a3})m\]
\[K_{11}(m) = \imath(-1)^{deg(m)}(\delta_{b0} - \delta_{b3})m\]
\[K_{22}(m) = \imath(-1)^{deg(m)}(\delta_{c0} - \delta_{c3})m\]
As a consequence a monomial of multidegree $(a,b,c)$ will have   $Kw$ weight
\[Kw(m) = \imath(-1)^{deg(m)}(\delta_{a0} - \delta_{a3},\delta_{b0} - \delta_{b3},\delta_{c0} - \delta_{c3})\]
\qed
\begin{pro}
\label{pro:kadjoint}  The algebra $\ggr$ contains the $Kw$-homogeneous compontens of the operators $\imath L_j,\imath\Lambda_k$ and the $\imath \Lambda_k$. The decomposition of the $\imath L_j$ into $Kw$-homogeneous components is as follows: 
\[\imath L_0 = \imath L_0^{(0,0,0)} + \imath L_0^{(0,-\imath,0)} + \imath L_0^{(0,0,-\imath)} + \imath L_0^{(0,-\imath,-\imath)}\]
\[\imath L_1 = \imath L_1^{(0,0,0)} + \imath L_1^{(-\imath,0,0)} + \imath L_1^{(0,0,-\imath)} + \imath L_1^{(-\imath,0,-\imath)}\]
\[\imath L_2 = \imath L_2^{(0,0,0)} + \imath L_2^{(-\imath,0,0)} + \imath L_2^{(0,-\imath,0)} + \imath L_0^{(-\imath,-\imath,0)}\]
The corresponding decomposition for the $\imath\Lambda_j$ has the same form, with opposite eigenvalues. 
\end{pro}
\dimo 
The previous proposition allows us to  identify the components $Kw$-homogeneous of any given operator which is  multidegree homogeneous. In particular this applies to the $\imath L_j$, which turn out to have only nonzero $Kw$-homogeneous components with eigenvalues $0,-\imath$ and have therefore the  decomposition as in the statement of the proposition. 

The $\imath\Lambda_j$ have opposite multidegree and therefore opposite $Kw$ degree with respect to the $\imath L_j$.

If we indicate with $p(x)$ the real polynomial $x(x^2 +1)$, one has that the operators $p(ad(K_{mm}))$ vanish on the linear span of the $\imath L_j$. From the decomposition $1 = (x^2 +1) - x^2$ one obtains that for all $j,m$
\[\left(ad(K_{mm})^2 + Id\right)(\imath L_j) - ad(K_{mm})(\imath L_j) =\imath L_j\]
and moreover the first summand must be in the Kernel of $ad(K_{mm})$ and the second must be in the Kernel of $ad(K_{mm})^2 + 1$. The component of $\imath L_j$ in the Kernel of $ad(K_{mm})^2 + 1$ is the sum of an eigenvector for $ad(K_{mm})$ with eigenvalue $\imath$ and one with eigenvalue $-\imath$. The first of these two however was shown to be zero before. It follows that iterating this procedure for the subsequent values $0,1,2$ of $m$ one obtains all the possibly  nonzero   $Kw$-homogeneous components of the $\imath L_j$, which therefore lie inside $\ggr$. The same reasoning applies to the $\imath\Lambda_j$, whose $Kw$-homogeneous components are therefore in $\ggr$. 
\qed\\

\section{Bases for spaces of highest weight vectors}
\label{sec:bases}
In this section we find geometrically miningful bases for the highest weight vector spaces $HW_k$. These bases will allow us to study the algebra $\ggr$ more in detail,  when in the next section we will  restrict the newly generated operators to the subrepresentations $HW_k$.\\
To begin with, we describe the space $HW_0$, which is made of highest weight vectors for $\ggr$. It turns out that it is possible to use the operators $\imath L_j$ applied to $1\in\bigwedge^*T*X_p$ and the $\imath\Lambda_k$ applied to the volume form $Vol(W_0\oplus W_1\oplus W_2) = Vol(W_0)\wedge Vol(W_1)\wedge Vol(W_2)$ to build a canonical of $HW_0$.
\begin{pro}
\label{pro:basishw0} 
The elements of the form $p(\imath L_0,\imath L_1,\imath L_2)\cdot 1$, where $p$ is a monomial in three variables of degree less than or equal to $3$, form a basis for $HW_0^{ev}$. The elements of the form $p(\imath \Lambda_0,\imath \Lambda_1,\imath \Lambda_2)\cdot Vol(W_0\oplus W_1\oplus W_2)$, where $p$ is a monomial in three variables of degree less than or equal to $3$, form a basis for  $HW_0^{odd}$, corresponding to the even one by way of the Hodge $*$ operator.
\end{pro}
\dimo The dimension of $HW_0^{ev}$ (and $HW_0^{odd}$) is $20$ (from table \ref{tab:isotyp}), and the Hodge $*$ sends an element of the form  $p(\imath L_0,\imath L_1,\imath L_2)\cdot 1$ precisely to $p(\imath \Lambda_0,\imath \Lambda_1,\imath \Lambda_2)\cdot Vol(W_0\oplus W_1\oplus W_2)$ (from the proof of Theorem \ref{teo:superantihermitean}, which shows that $\imath \Lambda_j = * ~ \imath L_j ~ *$ and the fact that $* 1 = Vol(W_0\oplus W_1\oplus W_2)$). It is therefore enough to prove that the $20$ elements of the form $p(\imath L_0,\imath L_1,\imath L_2)\cdot 1$ obtained by varying $p$ among the monomials of degree less than or equal to three are all independent. To do that, we first observe that all these elements have different multidegree (as the $\imath L_j$ are multidegree homogeneous, with independent multidegree vectors). To conclude it is therefore enough to show that all the elements above are different from zero. Of course it is enought to check the non vanishing when $p$ has degree three, and using the symmetric group $\mathcal{S}_3$ we are reduced to prove that $-\imath L_0^2 L_1\cdot 1\not= 0$, $-\imath L_0^3\cdot 1\not= 0$ and $-\imath L_0 L_1 L_2\cdot 1 \not= 0$. As the above are all wedge operators, this is equivalent to showing that the following three forms are different from zero (at all points):
\[\omega_D\wedge\omega_D\wedge\omega_2,\quad \omega_D\wedge\omega_D\wedge\omega_D,\quad \omega_D\wedge\omega_1\wedge\omega_2\]
This can be easily checked using an adapted coframe like the one given by the $\{v_{ij}\}$, or one can look at Theorem 1.6 of \cite{G1}.
\qed\\

To describe a basis for $HW_1$ we can use the same reasoning as before, applying operators of the form  $p(\imath L_0,\imath L_1,\imath L_2)$ to the three forms $w_{10},w_{11},w_{12}$ defined in the proof of Proposition \ref{pro:table}, which are a basis for the space $HW_1\cap \bigwedge_\C^1T^*_pX$.
\begin{pro}
\label{pro:basishw1}
The 12 elements
\[ w_{10},\  \imath L_0 w_{10},\   \imath L_1 w_{10},\ \imath L_2 w_{10},\  L_0^2 w_{10},  \ L_1L_0 w_{10},\  L_2L_0 w_{10}\] \[ L_1^2 w_{10},\;  L_2L_1 w_{10} ,\; L_2^2 w_{10}, \; \imath L_0^3 w_{10}, \; \imath L_2L_1L_0 w_{10}   \]
together  with the other 24 elements in their \({\mathcal S}_3\) orbit, 
constitute a basis for $HW_1^{odd}$.  The image of this basis   under the Hodge $*$ constitutes a basis for $HW_1^{ev}$.
\end{pro}
\dimo
A simple computation in coordinates shows that with respect to the $\mathbf{sl}(4,\R)$ action the $(\mathbf{h}_1, \mathbf{h}_2, \mathbf{h}_2)$-weight of $w_{10}$ is $(-1,0,-2)$, it is a lowest weight vector and therefore the $\mathbf{sl}(4,\R)$ module that it generates is irreducible of dimension $36$, and hence its complexification  plus its image under the Hodge $*$ coincides with $HW_1$. Moreover, out of the ordered linear basis 
{\small \[\Lambda_0,\Lambda_1,\Lambda_2, [L_1,\Lambda_0],[L_1,\Lambda_2],[L_2,\Lambda_0],[L_2,\Lambda_1],H_0,H_1,H_2,[L_0,\Lambda_2],[L_0,\Lambda_1],L_0,L_1,L_2\]}
only the last $8$ elements do not annichilate $w_{10}$. Using the Poincar\'e-Birkhoff-Witt theorem to write a linear basis for $\mathcal{U}(\mathbf{sl}(4,\R))$, one sees that a set of generators for this module is made out of polynomials in the $L_j$ applied to the vectors $w_{10},w_{11}$ and $w_{12}$.
By considering the  known $(\mathbf{h}_1, \mathbf{h}_2, \mathbf{h}_2)$-weight decomposition of this $\mathbf{sl}(4,\R)$ module, and the fact that the operators $L_j$ have independent $(\mathbf{h}_1, \mathbf{h}_2, \mathbf{h}_2)$-weight, one finds out tha the only thig to check to prove linear independence of the elements listed in the the proposition is that the following elements (and hence also the elements in their $\mathcal{S}_3$ orbits)  are different from zero:
\[\imath L_0^3 w_{10},\quad \imath L_0L_1L_2 w_{10},\quad L_1^2 w_{10}\]
The non vanishing of the first one is immediate, as $\omega_D^3$ is a volume form for the distribution $W_1\oplus W_2$. The third one can be shown to be nonzero by considering that the $H_1$-weight of $w_{10}$ is $-2$ and therefore the  $\mathbf{sl}(2,\C)\cong <L_1,\Lambda_1,H_1>$ module that it generates must have dimension $3$. The last non-vanishing is a simple computation in coordinates (or can be similarly proved using a repeated application of $\mathbf{sl}(2,\C)$ actions).
\qed\\

To describe a basis for $HW_2$ we can again apply the same reasoning as before, acting with  operators of the form  $p(\imath L_0,\imath L_1,\imath L_2)$ on  the three forms $w_{10}\wedge w_{11},w_{10}\wedge w_{12}$ and $w_{11}\wedge w_{12}$, which are a basis for the space $HW_2\cap \bigwedge_\C^2T^*_pX$.
\begin{pro}
\label{pro:basishw2}
The $20$ elements
\[\begin{array}{lllllllll} 
w_{10}\wedge w_{11}  &\imath  L_0 w_{10}\wedge w_{11} &\imath  L_1 w_{10}\wedge w_{11} & 
\imath L_2 w_{10}\wedge w_{11}\\
L_0^2 w_{10}\wedge w_{11} & L_1L_0 w_{10}\wedge w_{11} & L_1^2 w_{10}\wedge w_{11}&  \\
w_{10}\wedge w_{12}  &\imath  L_0 w_{10}\wedge w_{12} &\imath  L_1 w_{10}\wedge w_{12} & 
\imath L_2 w_{10}\wedge w_{12}   \\
L_0^2 w_{10}\wedge w_{12}  & L_2^2 w_{10}\wedge w_{12} & & \\
w_{11}\wedge w_{12}  & \imath L_0 w_{11}\wedge w_{12} & \imath L_1 w_{11}\wedge w_{12} & 
\imath L_2 w_{11}\wedge w_{12}   \\
L_1^2 w_{11}\wedge w_{12} & L_1L_2 w_{11}\wedge w_{12} & L_2^2 w_{11}\wedge w_{12} &  \\
\end{array}\]
constitute a basis for $HW_2^{ev}$. They together with their image under the Hodge $*$ constitute a basis for $HW_2$.
\end{pro}
\dimo
A simple computation in coordinates shows that $w_{10}\wedge w_{11}$ is  a lowest weight vector of   $(\mathbf{h}_1, \mathbf{h}_2, \mathbf{h}_2)$-weight equal to $(0,-1,-1)$ with respect to the $\mathbf{sl}(4,\R)$ action, and therefore the $\mathbf{sl}(4,\R)$ module that it generates is irreducible of dimension $20$, and hence its complexification  plus its image under the Hodge $*$ coincides with $HW_2$.
Moreover, out of the ordered linear basis 
{\small \[\Lambda_0,\Lambda_1,\Lambda_2, [L_1,\Lambda_0],[L_1,\Lambda_2],[L_2,\Lambda_0],[L_2,\Lambda_1],H_0,H_1,H_2,[L_0,\Lambda_2],[L_0,\Lambda_1],L_0,L_1,L_2\]}
only the last $8$ elements do not annichilate $w_{10}\wedge w_{11}$. Using the Poincar\'e-Birkhoff-Witt theorem to write a linear basis for $\mathcal{U}(\mathbf{sl}(4,\R))$, one sees that a set of generators for this module is made out of polynomials in the $L_j$ applied to the vectors $w_{10}\wedge w_{11},w_{10}\wedge w_{12}$ and $w_{11}\wedge w_{12}$.\\
The non vanishing of $L_0^2w_{12}$ was proven in the previous propsition, and from this one gets immediately the non-vanishing  of $L_0^2w_{10}\wedge w_{12}$ and of all the other $5$  elements in its $\mathcal{S}_3$-orbit. Looking at the multidegree, one sees that all these nonzero elements are actually independent.\\
Among the remaining elements associated to polynomials homogeneous of degree $2$, there are $L_0L_1(w_{10}\wedge w_{11})$, $L_1L_2(w_{11}\wedge w_{12})$ and $L_2L_0(w_{12}\wedge w_{10})$, which form an $\mathcal{S}_3$-orbit and have all the same multidegree $(2,2,2)$.
We notice that there is the following relation (verifiable by hand in coordinates) among them:
\begin{equation}\label{eqn:relation}L_0L_1(w_{10}\wedge w_{11}) + L_1L_2(w_{11}\wedge w_{12}) + L_2L_0(w_{12}\wedge w_{10}) = 0\end{equation}
The same computation also shows as a byproduct that any two of them are independent, and therefore the space that the three of them span has exactly dimension two. Taking into account also the three generators (corresponding to polynomials of degree zero) we have up to now $8 +3 = 11$ independent vectors in $HW_2^{ev}$. To finish the proof it is enough to show that the nine elements corresponding to polynomials of degree one are all independent, as $HW_2^{ev}$ has dimension $20$. Using multigree, the action of $\mathcal{S}_3$  and the non-vanishing information on the elements of degree two, we are left with the following two verifications:
\[\imath  L_2(w_{10}\wedge w_{11}) \not= 0,\quad dim\left(<\imath  L_0(w_{10}\wedge w_{11}),\imath L_2(w_{10}\wedge w_{12})>\right) = 2 \]
Both this verifications are straightforward, as the forms involved have degree $4$ The first one can also be verified observing that the Lie algebra  $\mathbf{sl}(2,\C)\cong <L_2,\Lambda_2,H_2>$ acts on $w_{10}\wedge w_{11}$ generating an irreducible module of dimension $2$, as this element has $H_2$-weight equal to $-1$.
\qed\\

It remains to  be described a basis for $HW_3$. We can again apply the same reasoning as before, acting with  operators of the form  $p(\imath L_0,\imath L_1,\imath L_2)$ on  the  form $w_{10}\wedge w_{11}\wedge w_{12}$, which generates $HW_3\cap \bigwedge_\C^3T^*_pX$.
\begin{pro}
\label{pro:basishw3}
The elements
{\small \[\
w_{10}\wedge w_{11}\wedge w_{12} , \quad \;
\imath L_0(w_{10}\wedge w_{11}\wedge w_{12}), \quad  \; \imath  L_1(w_{10}\wedge w_{11}\wedge w_{12})
, \quad \;  \imath L_2(w_{10}\wedge w_{11}\wedge w_{12})
\]}
\noindent constitute a basis for $HW_3^{odd}$. They together with their image under the Hodge $*$ constitute a basis for $HW_3$.
\end{pro}
\dimo
A simple computation in coordinates shows that $w_{10}\wedge w_{11}\wedge w_{12}$ is  a lowest weight vector of   $(\mathbf{h}_1, \mathbf{h}_2, \mathbf{h}_2)$-weight equal to $(0,0,-1)$ with respect to the $\mathbf{sl}(4,\R)$ action, and therefore the $\mathbf{sl}(4,\R)$ module that it generates is irreducible of dimension $4$, and hence its complexification  plus its image under the Hodge $*$ coincides with $HW_3$. To show that the four elements listed are independent, it is enough to show that they are different from zero, as they have all different multidegrees. Using the action of $\mathcal{S}^3$ it is enough to show that $\imath L_2(w_{10}\wedge w_{11}\wedge w_{12})\not= 0$.
This follows by observing that the Lie algebra $\mathbf{sl}(2,\C)\cong <L_2,\Lambda_2,H_2>$ acts on $w_{10}\wedge w_{11}\wedge w_{12}$ generating an irreducible module of dimension $2$, as this element has $H_2$-weight equal to $-1$.
\qed\\

\section{Restriction of $\ggr$ to the isotypical components}
\label{sec:restriction}
Using the bases consctructed in the previous section, we will now provide an explicit description for the action of $\ggr$ on the spaces $HW_k$.
\begin{teo}
\label{teo:hw3}For every point $p\in X$,
the operators $V_0,V_1,V_2,A_0,A_1,A_2$ are zero when restricted to the complex vector  space $HW_3$ of $\mathbf{so}(3,\C)$ highest weight vectors inside the isotypical component of type $\rho_3$ of $\bigwedge^*_\C T^*_p X$. Moreover,
in the basis of Proposition \ref{pro:basishw3} one has that the operators $\imath L_0,\imath L_1,\imath L_2$ acting on $HW_3^{odd}$ have matrices:
{\small \[M_{\beta}(\imath L_0)=\begin{pmatrix}
0 & 0 & 0 &0\\
1 & 0 & 0 & 0\\
0 & 0 & 0 &0\\
0 & 0 & 0 &0\end{pmatrix},~
M_{\beta}(\imath L_1)=\begin{pmatrix}
0 & 0 & 0 &0\\
0 & 0 & 0 &0\\
1 & 0 & 0 & 0\\
0 & 0 & 0 &0\end{pmatrix},~
M_{\beta}(\imath L_2)=\begin{pmatrix}
0 & 0 & 0 &0\\
0 & 0 & 0 &0\\
0 & 0 & 0 &0\\
1 & 0 & 0 & 0\end{pmatrix}\]}
The restriction of $\ggr$ to $HW_3$ is isomorphic $\mathbf{sl}(4,\R)$. 
\end{teo}
\dimo
The vanishing of the operators $V_0,V_1,V_2,A_0,A_1,A_2$ on the basis elements is immediate for reasons of multidegree. The matrices for $\imath L_0,\imath L_1,\imath L_2$ in terms of the given basis are straightforward to compute. \\
The matrices of the $\imath \Lambda_j$  acting on $HW_3^{odd}$ are the transpose of  those of the corresponding $\imath L_j$  up to a nonzero scalar because $\imath \Lambda_j = -\left(\imath L_j\right)^\ss$ for $j =0,1,2$ and the basis is orthogonal with respect to the inner product for reasons of multidegree.
The restriction of $\ggr$ to $HW_3^{odd}$ is therefore isomorphic $\mathbf{sl}(4,\R)$.\\
From Proposition \ref{pro:sync} one has that also the restriction of $\ggr$ to $HW_3$ is isomorphic to  $\mathbf{sl}(4,\R)$ (as the map sending an operator $\phi$ to $-\phi^\ss$ is a Lie algebra isomorphism).
\qed
\begin{lem}
\label{lem:hwj}
For all $p\in X$ and for $j\in\{0,1,2\}$ the restriction of the Lie  algebra $\ggr^{ev}$ to the complex vector  space $HW_j^{ev}$ of even $\mathbf{so}(3,\C)$ highest weight vectors inside the isotypical component of type $\rho_j$ of $\bigwedge^*_\C T^*_p X$ includes the full complex special linear algebra: specifically it includes   $\mathbf{sl}(20,\C)$ when $j = 0$ (respectively   $\mathbf{sl}(36,\C)$ when $j=1$ and  $\mathbf{sl}(20,\C)$ when $j = 2$).
\end{lem}
\dimo In the following proof we obtain an explicit presentation of the restrictions of the algebra $\ggr^{ev}$ to the $HW_j^{ev}$. This is done making a ``root search" using as starting point the $Kw$-homogeneous components of the restrictions of the generators  $\imath L_0,\imath L_1,\imath L_2$, which lie in $\ggr$ from Proposition \ref{pro:kadjoint}. We list in the Appendix (Tables \ref{tab:0}, \ref{tab:2},  \ref{tab:1})  the matrices of these components, which turn out to be simple, having a few nonzero entries. 

To build the tables we proceded as follows: given a $Kw$-homoegebeous compontent $\imath L_j^{(u,v,z)}$ of an operator and a basis element $v$, we know that the multidegree of the form $\imath L_j^{(u,v,z)}(w)$ would be the sum of the multidegrees of $\imath L_j^{(u,v,z)}$ and of $v$. Using this information and  Proposition \ref{pro:kdegree} we obtain vanishing of many candidate entries. When the $Kw$-degrees match, we check whether $\imath L_j(v)$ is another basis element: in this case we get non-vanishing. This is sufficient, except in a specific case in $HW_2$ when we have to  use  relation \ref{eqn:relation} from the proof of Proposition \ref{pro:basishw2} to express $\imath L_j^{(u,v,z)}(w)$ in terms of the basis. In the other cases the non-vanishing is verified directly, but this verification   turns out to be irrelevant  for the root search; this fact was signalled by putting the symbol $\times$ in the table (see the Appendix). \\
The completed root search proves that the restriction of $\ggr$ to $HW_j$ (for $j = 0,1,2$) includes a split real form for the full complex special linear algebra. As the Cartan contains the purely imaginary elements $K_{ii}$, all the weight spaces which have  $Kw$-weight different from $(0,0,0)$ must have real dimension two, and therefore the restriction must include the full complex special linear algebra.
\qed
\begin{teo}
\label{teo:hw2}
For all $p\in X$ and for $j\in\{0,1,2\}$ the restriction of the Lie superalgebra algebra $\ggr$ to the complex vector  space $HW_j$ of $\mathbf{so}(3,\C)$ highest weight vectors inside the isotypical component of type $\rho_j$ of $\bigwedge^*_\C T^*_p X$ is the Lie superalgebra $\mathbf{su}(HW_j^{ev}|HW_j^{odd},<~,~>)$ of all the operators which preserve the odd nondegenerate super Hermitean inner product $<~,~>_p$.
\end{teo}
\dimo  
From Propositions \ref{pro:injectively} and \ref{pro:sync}, we know how the action of $\ggr^{ev}$ on the $HW^{odd}_j$ is determined by the action on $HW^{ev}_j$: in particular from the previous lemma we know that $\ggr^{ev}$ is mapped isomorphically onto the even part of $\mathbf{su}(HW_j^{ev}|HW_j^{odd},<~,~>)$. To conclude the proof, it remains to be proven that $\ggr|_{HW_j}$ contains all the odd part of $\mathbf{su}(HW_j^{ev}|HW_j^{odd},<~,~>)$ and also the elements
\[\begin{pmatrix} \imath I & 0\\ 0 & \imath I\end{pmatrix}\]
Once obtained all the odd elements of $\mathbf{su}(HW_j^{ev}|HW_j^{odd},<~,~>)$, the one above can be obtained simply via the superbracket
\[\begin{pmatrix} \imath I & 0\\ 0 & \imath I\end{pmatrix} = \frac{1}{2}\left[\begin{pmatrix} 0 &  I & \\ \imath I & 0\end{pmatrix},\begin{pmatrix} 0 &  I & \\ \imath I & 0\end{pmatrix}\right]\]
We first observe that the odd operators of  $\mathbf{su}(HW_j^{ev}|HW_j^{odd},<~,~>)$ which send even forms into odd ones and the operators which send odd forms into even ones are two  modules of the same dimension for the adjoint action of the Lie algebra $\mathbf{su}(HW_j^{ev}|HW_j^{odd},<~,~>)^{ev}$ on $\mathbf{su}(HW_j^{ev}|HW_j^{odd},<~,~>)^{odd}$, which for simplicity we indicate in this proof respectively with   $An_j$ and $He_j$ (as, from Proposition \ref{pro:sync}, they are made respectively of matrices which have anti-hermitean and hermitean blocks with respect to any real orthonormal basis). These two modules have a decomposition (as real vector spaces and also as $\mathbf{sl}(n_j,\R)$-modules, where $n_0 = 20,~n_1 = 36$ and $n_2 = 20$) 
\[He_j = Re(He_j)\oplus_\R \imath Im(He_j),\quad An_j = Re(An_j)\oplus_\R \imath Im(An_j)\]
Given an element of $\ggr^{ev}|_{HW_j}$ with components $\phi$ and $-\phi^*$ which act on the even and on the odd parts of $HW_j$, its adjoint action on the odd part of the algebra is as follows:
\[ \left[\begin{pmatrix} \phi & 0 \\0 & -\phi^*\end{pmatrix},\begin{pmatrix}0 & B \\ A & 0\end{pmatrix}\right] = \begin{pmatrix} 0 & \phi B + B\phi^*\\ -\phi^*A - A\phi & 0\end{pmatrix}\]
from which it follows easily that considering the action of the Lie subalgebra $\mathbf{sl}(n_j,\R)\subset \ggr^{ev}|_{HW_j}$  given by real operators $\phi$, the four modules indicated above are all ireducible, of dimensions respectively $n_j(n_j +1)/2$, $n_j(n_j - 1)/2$,  $n_j(n_j - 1)/2$, $n_j(n_j + 1)/2$. As all the $\imath V_j$ have nonzero components exactly in $\imath Im(He_j)$ and in $\imath Im(An_j)$ simply by taking into account their multidegree (and the fact that they are purely imaginary), by acting on them with the algebra $\mathbf{sl}(n_j,\R)$ we can generate all $\imath Im(He_j)$ and $\imath Im(An_j)$ (which are not isomorphic having different dimension). Similarly, by acting with the algebra $\mathbf{sl}(n_j,\R)$ on the operators $A_j$ we generate all $Re(He_j)$ and $Re(An_j)$.
\qed

\begin{teo}
\label{teo:mainteo} The restrictions to the $HW_j$ provide an isomorphism of  Lie superalgebras
\[ 
\ggr~=~\mathbf{su}(HW_0^{ev}|HW_0^{odd},<~,~>_p)~ \oplus~  \mathbf{su}(HW_1^{ev}|HW_1^{odd},<~,~>_p)~ \oplus\]
\[ 
\oplus ~ 
\mathbf{su}(HW_2^{ev}|HW_2^{odd},<~,~>_p)~ \oplus~ \mathbf{su}(HW_3^{ev}|HW_3^{odd},<~,~>_p,\R)^{ev}\cong\]
\[\cong 
\mathbf{su}(20|20,<~,~>)~\oplus~ \mathbf{su}(36|36,<~,~>)~\oplus~ \mathbf{su}(20|20,<~,~>)~\oplus~ \mathbf{sl}(4,\R)\]
\end{teo}
\dimo The only thing remaining to be proven is that the restriction map of $\ggr$ to on $HW_0\oplus HW_2$ is surjective onto 
\[\mathbf{su}(HW_0^{ev}|HW_0^{odd},<~,~>)\oplus \mathbf{su}(HW_2^{ev}|HW_2^{odd},<~,~>)\] To do this, it is enough to find an even operator which restricts to zero on one module and is nonzero once restricted to the other. Ad example of these is given by the $K_{ij}$, which are nonzero on $HW_0$ and are zero on $HW_2$ for reasons of multidegree.
\qed
\begin{rmk} The (real) dimension of $\ggr$ is $8396$. This makes precise the observations at the end of Section \ref{sec:so2R}.
\end{rmk}
\begin{cor} The complexification $\gg$ of $\ggr$ is isomorphic via the restriction maps to the complex Lie superalgebra 
\[\mathbf{sl}(20|20)~\oplus~ \mathbf{sl}(36|36)~\oplus~ \mathbf{sl}(20|20)~\oplus~ \mathbf{sl}(4,\C)\]
\end{cor}
\section{Final remarks}
\label{sec:finalremarks}
We have at this point enough information to prove that our complexified algebra $\gg$ is indeed a $*$-Lie superalgebra, with respect to the standard adjunction operator $\dag$ associated to the natural superHilbert space structure on ${\bigwedge}_\C^*T_p^*X$. For conveninece of the reader we reproduce here the definition of superHilbert space and $*$ operator (induced by the superadjunction $\dag$), taken from  \cite{CCTV} and \cite{V}:
\begin{dfn}[\cite{CCTV}]
A {\em super Hilbert space}
 is a super vector space
${\hh}={\hh}_0\oplus {\hh}_1$ over $\C$ with a
scalar product $(~ ,~)$ such that
${\hh}$ is a Hilbert space
under $(~ ,~)$, and ${\hh}_i (i=0,1)$
are mutually orthogonal closed linear subspaces.
If we define
$$
\langle x, y\rangle=
\begin{cases}
0 & \mbox{ if } x \mbox{ and } y \mbox{ are of opposite parity}\cr
(x, y) & \mbox{ if } x \mbox{ and } y \mbox{ are even } \cr
i(x,y) & \mbox{ if } x \mbox{ and } y \mbox{ are odd }  \cr
\end{cases}
$$
then $\langle x,y\rangle$ is an even super Hermitean form with
$$
\langle y,x\rangle=
(-1)^{deg(x)deg(y)}\overline {\langle x,y\rangle},\
\langle x,x\rangle >0 ~(x\not=0\ {\rm  even }),\
i^{-1}\langle x,x\rangle >0
~(x\not=0\  {\rm  odd }).
$$
If $T({\hh}\rightarrow {\hh})$ is a bounded linear operator,
we denote by $T^\ast$ its Hilbert
space adjoint
and by $T^{\dag} $ its super adjoint given by $\langle Tx,y\rangle =
(-1)^{deg(T)deg(x)}\langle x, T^{\dag} y\rangle $.
\end{dfn}
We remark that our spaces of sections form a preHilbert space, and would need to be completed to become a Hilbert space proper.\\
We proved in Theorem \ref{teo:mainteo} that  our real Lie superalgebra $\ggr$ is the full real Lie superalgebra  of all the odd-unitary operators with respect to the odd super Hermitean inner product $<~,~>$ for all the $\mathbf{so}(3,\R)$-isotypical components, save for the smallest one where we obtain only the real and even part of it. From this it follows that to prove that $\ggr$ is closed with respect to $\dag$ it is enough to show that an element of the form 
\[\begin{pmatrix} A & B =~ ^t{\overline{B}}\\ C =~ -{^t{\overline{C}}} & -{^t{\overline{A}}}\end{pmatrix}\]
is sent to an element of the same form by $\dag$. Notice that we are working on a fixed isotypical component, and this is allowed because the different components are orthogonal with respect to the inner product $\langle ~,~\rangle$.
It is now easy to compute:
\[\begin{pmatrix} A & B =~ ^t{\overline{B}}\\ C =~ -{^t{\overline{C}}} & -{^t{\overline{A}}}\end{pmatrix}^\dag = \begin{pmatrix} ^t\overline{A} & \imath C \\ -\imath B & -A\end{pmatrix}\]
and therefore if $\phi\in\ggr$ then $\phi^\dag \in\ggr$. This implies that also $\gg$ is closed with respect to $\dag$,  as was to be proven.

\[\]
\[\]

\appendix
\section{Tables}
The notation in the tables is explained by the following example:  an entry of the form $L_0^{(0,0,-)}$ indicates that the matrix of the $Kw$-homogeneous component $\imath L_0^{(0,0,-\imath)}$ of the operator $\imath L_0$ of $Kw$-degree $(0,0,-\imath)$ has a  nonzero entry in that position. The basis elements have been indicated as follows: $(a,b,c)$ indicates the form $(\imath)^{a+b+c}L_0^aL_1^bL_2^c \ 1$ (Table \ref{tab:0}), $(a,b,c)j$ indicates the form $(\imath)^{a+b+c}L_0^aL_1^bL_2^c \  w_{1j}$ (Tables \ref{tab:1}, \ref{tab:1-2}) and  $(a,b,c)jk$ indicates the form $(\imath)^{a+b+c}L_0^aL_1^bL_2^c \ w_{1j}\wedge w_{1k}$ (Table \ref{tab:2}).\\
When an entry is followed by the symbol $\times$ it means that the component is present with a nonzero coefficient, that this fact was proven by a direct ``brute force" computation and that it is in any case irrelevant to know whether this coefficient is zero or not.\\
\newpage
\pagestyle{empty}

{\tiny
\begin{sidewaystable}
\centering
\begin{tabular}{|p{0.9cm}|p{0.9cm}|p{0.9cm}|p{0.9cm}|p{0.9cm}|p{0.9cm}|p{0.9cm}|p{0.9cm}|p{0.9cm}|p{0.9cm}|p{0.9cm}|p{0.9cm}|p{0.8cm}|p{0.8cm}|p{0.8cm}|p{0.8cm}|p{0.8cm}|p{0.8cm}|p{0.8cm}|p{0.8cm}|p{0.8cm}|}
\hline
   &(0,0,0)  & (1,0,0) &(0,1,0) & (0,0,1)&(2,0,0) & (1,1,0)	& (1,0,1)&	(0,2,0)&	(0,1,1)&	(0,0,2)&	(3,0,0)	&(2,1,0)&	(2,0,1)	&(1,2,0)&	(1,1,1)&	(1,0,2)&	(0,3,0)&	(0,2,1)	&(0,1,2)	& (0,0,3) \\
 \hline
 (0,0,0) &  & & & & & & & & & & & & & & & & & & &	\\
 \hline
(1,0,0)&	$L_0^{(0,-,-)}$	&	& & & & & & & & & & & & & & & & & &	\\
 \hline
(0,1,0)&	$L_1^{(-,0,-)}$	& & & & & & & & & & & & & & & & & & &	\\
 \hline
(0,0,1)&	$L_2^{(-,-,0)}$		&	& & & & & & & & & & & & & & & & & &	\\
 \hline
(2,0,0)&		& $L_0^{(0,0,0)}$	&	& & & & & & & & & & & & & & & & &	\\
 \hline
(1,1,0)	& &	$L_1^{(-,0,0)}$	& $L_0^{(0,-,0)}$ & & & & & & & & & & & & & & & & &	\\
 \hline
(1,0,1)	& &	$L_2^{(-,0,0)}$	& & 	$L_0^{(0,0,-)}$	& & & & & & & & & & & & & & & &  \\
 \hline
 (0,2,0)	& & &		$L_1^{(0,0,0)}$	& & & & & & & & & & & & & & & & & \\
 \hline
(0,1,1)& & & 			$L_2^{(0,-,0)}$	& $L_1^{(0,0,-)}$	& 	& & & & & & & & & & & & & & &  \\
 \hline
(0,0,2)	& & & &			$L_2^{(0,0,0)}$	& & & & & & & & & & & & & & & &  \\
 \hline
(3,0,0)	& & & & &				$L_0^{(0,-,-)}$& 	& & & & & & & & & & & & & & \\
 \hline
 (2,1,0)	& & & & &				$L_1^{(-,0,-)}$	& $L_0^{(0,0,-)}$ & & & & & & & & & & & & & & \\
 \hline
(2,0,1)	& & & & &				$L_2^{(-,-,0)}$	& & 	$L_0^{(0,-,0)}$	& 	& & & & & & & & & & & &  \\
 \hline
 (1,2,0)	& & & & & &					$L_1^{(0,0,-)}$	&  &	$L_0^{(0,-,-)}$	& 	& & & & & & & & & & &  \\ 
 \hline
(1,1,1)		& & & & & &				$L_2^{(0,0,0)}$	& $L_1^{(0,0,0)}$	& &	$L_0^{(0,0,0)}$	& & & & & & & & & & &  \\
 \hline
 (1,0,2) & & & & & &  &							$L_2^{(0,-,0)}$	& & & 		$L_0^{(0,-,-)}$	& 	& & & & & & & & & \\
 \hline
 (0,3,0)	&  & & &  &  & &  &							$L_1^{(-,0,-)}$	& &  & 		$K_{01}$  & 	& & & & & & & & \\
 \hline
 (0,2,1)		&  & & & & & & &						$L_2^{(-,-,0)}$& 	$L_1^{(-,0,0)}$		& 	& & & & & & & & & &	 \\
 \hline
(0,1,2)		&  & & & & & & &	&							$L_2^{(-,0,0)}$& 	$L_1^{(-,0,-)}$	& 	& & & & & & & & &\\
 \hline
 (0,0,3)	& & & & & & & & & &									$L_2^{(-,-,0)}$	&$K_{02}$	& & & & & &				$K_{12}$	& & &  \\
\hline
\end{tabular}
\medskip
\caption{Table for $HW_0$}
\label{tab:0}
\end{sidewaystable}}

\newpage
\topmargin=0cm
\thispagestyle{empty}
{\tiny
\begin{sidewaystable}
\centering
\begin{tabular}{|p{0.95cm}|p{0.95cm}|p{0.95cm}|p{0.95cm}|p{0.95cm}|p{0.8cm}|p{0.8cm}|p{0.8cm}|p{0.95cm}|p{0.95cm}|p{0.95cm}|p{0.95cm}|p{0.8cm}|p{0.8cm}|p{0.95cm}|p{0.95cm}|p{0.95cm}|p{0.95cm}|p{0.8cm}|p{0.8cm}|p{0.8cm}|}
\hline
&	(0,0,0)01  &	(1,0,0)01	& (0,1,0)01	& (0,0,1)01	&(2,0,0)01	&(1,1,0)01&	(0,2,0)01&	(0,0,0)02&	(1,0,0)02&	(0,1,0)02&	(0,0,1)02&	(2,0,0)02	&(0,0,2)02&	(0,0,0)12	&(1,0,0)12	&(0,1,0)12	&(0,0,1)12&	(0,2,0)12&	(0,1,1)12&	(0,0,2)12 \\
 \hline
(0,0,0)01 & & & & &	 & & & & & & & & & &		& & & & &	\\
 \hline
(1,0,0)01& 	$L_0^{(0,0,-)}$	& & & & & & & & & & & & & & & & & & & \\
 \hline
(0,1,0)01 & 	$L_1^{(0,0,-)}$	& & & & & & & & & & & & & & & & & & &  \\
 \hline
(0,0,1)01 &	$L_2^{(0,0,0)}$	& & & & & & & & & & & & & & & & & & &  \\
 \hline
(2,0,0)01	&	&$L_0^{(0,-,0)}$	&	&	&	&	&	&	&	&	&	&	&	&	&$L_2^{(-,-,0)}$x	&	&$L_0^{(0,-,0)}$x	&	&	&	\\
 \hline
(1,1,0)01	&	&$L_1^{(0,0,0)}$	&$L_0^{(0,0,0)}$	&	&	&	&	&	&$L_2^{(0,0,0)}$	&	&$L_0^{(0,0,0)}$	&	&	&	&	&	&	&	&	&	\\
 \hline
(0,2,0)01	&	&	&$L_1^{(-,0,0)}$	&	&	&	&	&	&	& $L_2^{(-,-,0)}$x	&$L_1^{(-,0,0)}$x	&	&	&	&	&	&	&	&	&  \\
 \hline
(0,0,0)02	&	&	& &	&	&	&	&		&	&	&	&	&	&	&	&	&	&	&	& \\
 \hline
(1,0,0)02	&	&	&	&	&	&	&	&$L_0^{(0,-,0)}$	&	&	&	&	&	&	&	&	&	&	&	&  \\
 \hline
(0,1,0)02	&	&	&	&	&	&	&	&$L_1^{(0,0,0)}$	&	&	&	&	&	&	&	&	&	&	&	&  \\
 \hline
(0,0,1)02	&	&	&	&	&	&	&	&$L_2^{(0,-,0)}$	&	&	&	&	&	&	&	&	&	&	&	& \\
 \hline
(2,0,0)02	&	&	&	&	&	&	&	&	&$L_0^{(0,0,-)}$	&	&	&	&	&	&$L_1^{(-,0,-)}$x	&$L_0^{(0,0,-)}$x	&	&	&	& \\
 \hline
(0,0,2)02	&	&	&$L_2^{(-,0,0)}$x	&$L_1^{(-,0,-)}$x	&	&	&	&	&	&	&$L_2^{(-,0,0)}$	&	&	&	&	&	&	&	&	&  \\
 \hline
(0,0,0)12	&	&	&	&	&	&	&	&	&	&	&	&	&	&	&	&	&	&	&	& \\ 
 \hline
(1,0,0)12	&	&	&	&	&	&	&	&	&	&	&	&	&	&$L_0^{(0,0,0)}$	&	&	&	&	&	& \\
 \hline
(0,1,0)12	&	&	&	&	&	&	&	&	&	&	&	&	&	&$L_1^{(-,0,0)}$	&	&	&	&	&	& \\
 \hline
(0,0,1)12	&	&	&	&	&	&	&	&	&	&	&	&	&	&$L_2^{(-,0,0)}$	&	&	&	&	&	& \\
  \hline
(0,2,0)12	&	&	&	&	&	&	&	&	&$L_1^{(0,0,-)}$x	&$L_0^{(0,-,-)}$x	&	&	&	&	&	&$L_1^{(0,0,-)}$	&	&	&	&\\
 \hline
(0,1,1)12	&	&	&	&	&	&	&	&	&$L_2^{(0,0,0)}$	&	&$L_0^{(0,0,0)}$	&	&	&	&	&$L_2^{(0,0,0)}$	&$L_1^{(0,0,0)}$	&	&	&\\
\hline
(0,0,2)12	&	&$L_2^{(0,-,0)}$x	&	&$L_0^{(0,-,-)}$x	&	&	&	&	&	&	&	&	&	&	&	&	&$L_2^{(0,-,0)}$	&	&	&\\
\hline
\end{tabular}
\medskip
\caption{Table for $HW_2$}
\label{tab:2}
\end{sidewaystable}}

\newpage
\thispagestyle{empty}

{\tiny
\begin{sidewaystable}
\centering
\begin{tabular}{|p{0.95cm}|p{0.95cm}|p{0.95cm}|p{0.95cm}|p{0.95cm}|p{1.0cm}|p{1.0cm}|p{1.0cm}|p{1.0cm}|p{0.95cm}|p{1.0cm}|p{0.85cm}|p{0.85cm}|p{0.95cm}|p{0.95cm}|p{0.95cm}|p{1cm}|p{1cm}|p{0.95cm}|}
\hline
&	(0,0,0)0	&	(1,0,0)0	&	(0,1,0)0	&	(0,0,1)0	&	(2,0,0)0	&	(1,1,0)0&		(1,0,1)0&		(0,2,0)0&		(0,1,1)0&		(0,0,2)0&		(3,0,0)0	&	(1,1,1)0&		(0,0,0)1&		(1,0,0)1&		(0,1,0)1	&	(0,0,1)1	&	(2,0,0)1&	(1,1,0)1	\\ \hline
(0,0,0)0	&	&	&	&	&	&	&	&	&	&	&	&	&	&	&	&	&	&	\\	\hline
(1,0,0)0	&$L_0^{(0,-,-)}$	&	&	&	&	&	&	&	&	&	&	&	&	&	&	&	&	&	\\	\hline
(0,1,0)0	&$L_1^{(0,0,-)}$	&	&	&	&	&	&	&	&	&	&	&	&	&	&	&	&	&	\\	\hline
(0,0,1)0	&$L_2^{(0,-,0)}$	&	&	&	&	&	&	&	&	&	&	&	&	&	&	&	&	&	\\	\hline
(2,0,0)0	&	&$L_0^{(0,0,0)}$	&	&	&	&	&	&	&	&	&	&	&	&	&	&	&	&	\\	\hline
(1,1,0)0	&	&$L_1^{(0,0,0)}$	&$L_0^{(0,-,0)}$	&	&	&	&	&	&	&	&	&	&	&	&	&	&	&	\\	\hline
(1,0,1)0	&	&$L_2^{(0,0,0)}$	&	&$L_0^{(0,0,-)}$	&	&	&	&	&	&	&	&	&	&	&	&	&	&	\\	\hline
(0,2,0)0	&	&	&$L_1^{(-,0,0)}$	&	&	&	&	&	&	&	&	&	&	&	&	&	&$K_{01}$	&	\\	\hline
(0,1,1)0	&	&	&$L_2^{(-,-,0)}$	&$L_1^{(-,0,-)}$	&	&	&	&	&	&	&	&	&	&	&	&	&	&	\\	\hline
(0,0,2)0	&	&	&	&$L_2^{(-,0,0)}$	&	&	&	&	&	&	&	&	&	&	&	&	&	&	\\	\hline
(3,0,0)0	&	&	&	&	&$L_0^{(0,-,-)}$	&	&	&	&	&	&	&	&	&	&	&	&$L_1^{(-,0,-)}$x	&$L_0^{(0,-,-)}$x	\\	\hline
(1,1,1)0	&	&	&	&	&	&$L_2^{(-,0,0)}$	&$L_1^{(-,0,0)}$	&	&$L_0^{(0,0,0)}$	&	&	&	&	&	&	&	&	&	\\	\hline
(0,0,0)1	&	&	&	&	&	&	&	&	&	&	&	&	&	&	&	&	&	&	\\	\hline
(1,0,0)1	&	&	&	&	&	&	&	&	&	&	&	&	&$L_0^{(0,0,-)}$	&	&	&	&	&	\\	\hline
(0,1,0)1	&	&	&	&	&	&	&	&	&	&	&	&	&$L_1^{(-,0,-)}$	&	&	&	&	&	\\	\hline
(0,0,1)1	&	&	&	&	&	&	&	&	&	&	&	&	&$L_2^{(-,0,0)}$	&	&	&	&	&	\\	\hline
(2,0,0)1	&	&	&	&	&	&	&	&$K_{10}$&	&	&	&	&	&$L_0^{(0,-,0)}$	&	&	&	&	\\	\hline
(1,1,0)1	&	&	&	&	&	&	&	&	&	&	&	&	&	&$L_1^{(-,0,0)}$	&$L_0^{(0,0,0)}$	&	&	&	\\	\hline
(1,0,1)1	&	&	&	&	&	&	&	&	&	&	&	&	&	&$L_2^{(-,-,0)}$	&	&$L_0^{(0,-,-)}$	&	&	\\	\hline
(0,2,0)1	&	&	&	&	&	&	&	&	&	&	&	&	&	&	&$L_1^{(0,0,0)}$	&	&	&	\\	\hline
(0,1,1)1	&	&	&	&	&	&	&	&	&	&	&	&	&	&	&$L_2^{(0,0,0)}$	&$L_1^{(0,0,-)}$	&	&	\\	\hline
(0,0,2)1	&	&	&	&	&	&	&	&	&	&	&	&	&	&	&	&$L_2^{(0,-,0)}$	&	&	\\	\hline
(0,3,0)1	&	&	&	&	&	&$L_1^{(-,0,-)}$x	&	&$L_0^{(0,-,-)}$x	&	&	&	&	&	&	&	&	&	&	\\	\hline
(1,1,1)1	&	&	&	&	&$L_2^{(0,-,0)}$x	&	&$L_0^{(0,-,0)}$x	&	&	&	&	&	&	&	&	&	&	&$L_2^{(0,-,0)}$	\\	\hline
(0,0,0)2	&	&	&	&	&	&	&	&	&	&	&	&	&	&	&	&	&	&	\\	\hline
(1,0,0)2	&	&	&	&	&	&	&	&	&	&	&	&	&	&	&	&	&	&	\\	\hline
(0,1,0)2	&	&	&	&	&	&	&	&	&	&	&	&	&	&	&	&	&	&	\\	\hline
(0,0,1)2	&	&	&	&	&	&	&	&	&	&	&	&	&	&	&	&	&	&	\\	\hline
(2,0,0)2	&	&	&	&	&	&	&	&	&	&$K_{20}$	&	&	&	&	&	&	&	&	\\	\hline
(1,1,0)2	&	&	&	&	&	&	&	&	&	&	&	&	&	&	&	&	&	&	\\	\hline
(1,0,1)2	&	&	&	&	&	&	&	&	&	&	&	&	&	&	&	&	&	&	\\	\hline
(0,2,0)2	&	&	&	&	&	&	&	&	&	&	&	&	&	&	&	&	&	&	\\	\hline
(0,1,1)2	&	&	&	&	&	&	&	&	&	&	&	&	&	&	&	&	&	&	\\	\hline
(0,0,2)2	&	&	&	&	&	&	&	&	&	&	&	&	&	&	&	&	&	&	\\	\hline
(0,0,3)2	&	&	&	&	&	&	&$L_2^{(-,-,0)}$x	&	&	&$L_0^{(0,-,-)}$x	&	&	&	&	&	&	&	&	\\	\hline
(1,1,1)2	&	&	&	&	&$L_1^{(0,0,-)}$x	&$L_0^{(0,0,-)}$x	&	&	&	&	&	&	&	&	&	&	&	&$L_1^{(0,0,-)}$x	\\	\hline
\end{tabular}
\medskip
\caption{Table for $HW_1$, first part}
\label{tab:1}
\end{sidewaystable}}

\newpage
\thispagestyle{empty}

{\tiny
\begin{sidewaystable}

\centering
\begin{tabular}{|p{0.95cm}|p{0.95cm}|p{1cm}|p{1cm}|p{1cm}|p{0.85cm}|p{0.85cm}|p{0.95cm}|p{0.95cm}|p{0.95cm}|p{0.95cm}|p{1cm}|p{0.95cm}|p{1cm}|p{1cm}|p{1cm}|p{1cm}|p{0.85cm}|p{0.85cm}|}
	\hline
	&(1,0,1)1	&(0,2,0)1	&(0,1,1)1&	(0,0,2)1&	(0,3,0)1	&(1,1,1)1	&(0,0,0)2&	(1,0,0)2	&(0,1,0)2	&(0,0,1)2	&(2,0,0)2	&(1,1,0)2	&(1,0,1)2&	(0,2,0)2	&(0,1,1)2	&(0,0,2)2	&(0,0,3)2	&(1,1,1)2	\\	\hline
(0,0,0)0	&	&	&	&	&	&	&	&	&	&	&	&	&	&	&	&	&	&	\\	\hline
(1,0,0)0	&	&	&	&	&	&	&	&	&	&	&	&	&	&	&	&	&	&	\\	\hline
(0,1,0)0	&	&	&	&	&	&	&	&	&	&	&	&	&	&	&	&	&	&	\\	\hline
(0,0,1)0	&	&	&	&	&	&	&	&	&	&	&	&	&	&	&	&	&	&	\\	\hline
(2,0,0)0	&	&	&	&	&	&	&	&	&	&	&	&	&	&	&	&	&	&	\\	\hline
(1,1,0)0	&	&	&	&	&	&	&	&	&	&	&	&	&	&	&	&	&	&	\\	\hline
(1,0,1)0	&	&	&	&	&	&	&	&	&	&	&	&	&	&	&	&	&	&	\\	\hline
(0,2,0)0	&	&	&	&	&	&	&	&	&	&	&	&	&	&	&	&	&	&	\\	\hline
(0,1,1)0	&	&	&	&	&	&	&	&	&	&	&	&	&	&	&	&	&	&	\\	\hline
(0,0,2)0	&	&	&	&	&	&	&	&	&	&	&	&	&	&	&	&	&	&	\\	\hline
(3,0,0)0	&	&	&	&	&	&	&	&	&	&	&$L_2^{(-,-,0)}$x	&	&$L_0^{(0,-,-)}$x	&	&	&	&	&	\\	\hline
(1,1,1)0	&	&$L_2^{(-,0,0)}$x	&$L_1^{(-,0,0)}$x	&	&	&	&	&	&	&	&	&	&	&	&$L_2^{(-,0,0)}$x	&$L_1^{(-,0,0)}$x	&	&	\\	\hline
(0,0,0)1	&	&	&	&	&	&	&	&	&	&	&	&	&	&	&	&	&	&	\\	\hline
(1,0,0)1	&	&	&	&	&	&	&	&	&	&	&	&	&	&	&	&	&	&	\\	\hline
(0,1,0)1	&	&	&	&	&	&	&	&	&	&	&	&	&	&	&	&	&	&	\\	\hline
(0,0,1)1	&	&	&	&	&	&	&	&	&	&	&	&	&	&	&	&	&	&	\\	\hline
(2,0,0)1	&	&	&	&	&	&	&	&	&	&	&	&	&	&	&	&	&	&	\\	\hline
(1,1,0)1	&	&	&	&	&	&	&	&	&	&	&	&	&	&	&	&	&	&	\\	\hline
(1,0,1)1	&	&	&	&	&	&	&	&	&	&	&	&	&	&	&	&	&	&	\\	\hline
(0,2,0)1	&	&	&	&	&	&	&	&	&	&	&	&	&	&	&	&	&	&	\\	\hline
(0,1,1)1	&	&	&	&	&	&	&	&	&	&	&	&	&	&	&	&	&	&	\\	\hline
(0,0,2)1	&	&	&	&	&	&	&	&	&	&	&	&	&	&	&	&	&	&	\\	\hline
(0,3,0)1	&	&$L_1^{(-,0,-)}$	&	&	&	&	&	&	&	&	&	&	&	&$L_2^{(-,-,0)}$x	&$L_1^{(-,0,-)}$x	&	&	&	\\	\hline
(1,1,1)1	&$L_1^{(0,0,0)}$	&	&$L_0^{(0,-,0)}$	&	&	&	&	&	&	&	&	&	&$L_2^{(0,-,0)}$x	&	&	&$L_0^{(0,-,0)}$x	&	&	\\	\hline
(0,0,0)2	&	&	&	&	&	&	&	&	&	&	&	&	&	&	&	&	&	&	\\	\hline
(1,0,0)2	&	&	&	&	&	&	&$L_0^{(0,-,0)}$	&	&	&	&	&	&	&	&	&	&	&	\\	\hline
(0,1,0)2	&	&	&	&	&	&	&$L_1^{(-,0,0)}$	&	&	&	&	&	&	&	&	&	&	&	\\	\hline
(0,0,1)2	&	&	&	&	&	&	&$L_2^{(-,-,0)}$	&	&	&	&	&	&	&	&	&	&	&	\\	\hline
(2,0,0)2	&	&	&	&	&	&	&	&$L_0^{(0,0,-)}$	&	&	&	&	&	&	&	&	&	&	\\	\hline
(1,1,0)2	&	&	&	&	&	&	&	&$L_1^{(-,0,-)}$	&$L_0^{(0,-,-)}$	&	&	&	&	&	&	&	&	&	\\	\hline
(1,0,1)2	&	&	&	&	&	&	&	&$L_2^{(-,0,0)}$	&	&$L_0^{(0,0,0)}$	&	&	&	&	&	&	&	&	\\	\hline
(0,2,0)2	&	&	&	&$K_{21}$	&	&	&	&	&$L_1^{(0,0,-)}$	&	&	&	&	&	&	&	&	&	\\	\hline
(0,1,1)2	&	&	&	&	&	&	&	&	&$L_2^{(0,-,0)}$	&$L_1^{(0,0,0)}$	&	&	&	&	&	&	&	&	\\	\hline
(0,0,2)2	&	&	&	&	&	&	&	&	&	&$L_2^{(0,0,0)}$	&	&	&	&	&	&	&	&	\\	\hline
(0,0,3)2	&	&	&$L_2^{(-,-,0)}$x	&$L_1^{(-,0,-)}$x	&	&	&	&	&	&	&	&	&	&	&	&$L_2^{(-,-,0)}$	&	&	\\	\hline
(1,1,1)2	&	&$L_0^{(0,0,-)}$x	&	&	&	&	&	&	&	&	&	&$L_2^{(0,0,0)}$	&$L_1^{(0,0,-)}$	&	&$L_0^{(0,0,-)}$	&	&	&	\\	\hline
\end{tabular}
\medskip
\caption{Table for $HW_1$, second  part}
\label{tab:1-2}
\end{sidewaystable}}
\thispagestyle{empty}
\end{document}